\documentclass[11pt,a4paper,reqno]{amsart}

\usepackage[margin=3cm]{geometry}

\usepackage{epsfig}

\usepackage{xcolor}
\usepackage{psfrag}
\usepackage{comment}

\newtheorem{lm}{Lemma}
\newtheorem{theorem}{Theorem}

\newtheorem{remark}{Remark}
\def\fl{\!}

\def\ds{\displaystyle}

\def\qed{\hfill $\Box$ \\ \bigskip}

\newcommand{\Z}{\ensuremath{\mathbb{Z}}}

\newcommand{\sech}{\mathrm{sech}}

\newcommand{\Wuloc}{W^{\textrm{\small u}}_{\textrm{\small loc}}}
\newcommand{\Wsloc}{W^{\textrm{\small s}}_{\textrm{\small loc}}}

\definecolor{Sepia}{rgb}{0.6,0.3,0.1}
\definecolor{RawSienna}{rgb}{0.2,0.7,0.1}
\definecolor{NavyBlue}{rgb}{0.2,0,0.5}
\definecolor{PineGreen}{rgb}{0,0.5,0}
\definecolor{Sepia1}{rgb}{0.9,0.0,0.1}
\definecolor{Kirpich}{rgb}{0.8,0.1,0.1}
\definecolor{Sepia2}{rgb}{0.4,0.,0.4}
\definecolor{Sepia3}{rgb}{0.4,0.,0.}

\newcommand{\Mt}{\widetilde{M}}

\newcommand{\rme}{\mathrm{e}}
\newcommand{\rmi}{\mathrm{i}}

\newcommand{\Totkm}{T_{2m1k}}

\title[Attracting, repelling and elliptic orbits in reversible maps]
{Mixed dynamics of $2$-dimensional reversible maps with a symmetric couple of quadratic homoclinic tangencies}

\author{A. Delshams}
\address{Departament de Matem\`atiques, Universitat Polit\`ecnica de Catalunya, Barcelona, Spain}
\email{amadeu.delshams@upc.edu}

\author{M.S. Gonchenko}
\address{Departament de Matem\`atiques, Universitat de Barcelona, Barcelona, Spain}
\email{mgonchenko@gmail.com}

\author{S.V. Gonchenko}
\address{Lobachevsky University of Nizhny Novgorod,  N.Novgorod, Russia}
\email{sergey.gonchenko@mail.ru}

\author{J.T. L\'azaro}
\address{Departament de Matem\`atiques, Universitat Polit\`ecnica de Catalunya, Barcelona, Spain}
\email{jose.tomas.lazaro@upc.edu}

\thanks{
This work has been
supported by the Russian Scientific Foundation grant: sections 1--4, 6 and 7 were carried out under the project 14-41-00044, and section 5 under the project 14-12-00811.
AD, MG and JTL have been also partially supported
by the MICIIN/FEDER grant MTM2015-65715-P and by the Catalan grant 2014SGR-504 (AD, JTL).
MG has been partially supported by Juan de la Cierva-Formaci\'on Fellowship FJCI-2014-21229, the grant MTM2016-80117-P (MINECO/FEDER, UE) and the Knut 
and Alice Wallenberg Foundation grant 2013-0315. 
SG also thanks RFBR (grant 16-01-00364) and the Russian Ministry of Science and Education, project 1.3287.2017. }

\begin{document}

\begin{abstract}
We study dynamics and bifurcations of $2$-dimensional reversible
maps having a symmetric saddle fixed point with an asymmetric pair of nontransversal homoclinic orbits
(a symmetric nontransversal homoclinic figure-8).  We consider one-parameter
families of reversible maps unfolding the initial
homoclinic tangency and prove the existence of infinitely many
sequences (cascades) of bifurcations related to the birth of asymptotically
stable, unstable and elliptic periodic orbits.

\vskip1cm


\medskip

\end{abstract}

\maketitle

\section{Introduction}

The mathematical foundations of the Bifurcation Theory were laid in the  famous paper of Andronov and Pontryagin~\cite{AP37} where the notion of ``rough'' (structurally stable) systems was introduced. Later on,  in a series of (classical) papers by Andronov, Leontovich and Maier (see e.g. books~\cite{AGLM1,AGLM2}) it was proved that rough 2-dimensional systems form an open and dense set in the space of dynamical systems. This notion of roughness of a system (i.e. topological equivalence/conjugacy of the chosen system with any close system) is naturally extended to multidimensional systems. Such extension was carried out in the 60's
(this period was called by Anosov as the time of the ``hyperbolic revolution'') where structurally stable systems were also entitled as ``Hyperbolic Systems''. 

Such systems are divided in two large classes: Morse-Smale systems (with a simple dynamics) and hyperbolic systems with infinitely many periodic orbits. By definition, structural stable systems are open subsets. However, in the multidimensional case (that is, dimension
$\geq 3$ for flows and $\geq 2$ for diffeomorphisms), they are not dense, as it was first shown by Smale~\cite{Sm66,Sm67}.

A very important breakthrough was due to  Newhouse~\cite{N70,N79} who proved that, near any 2-dimensional diffeomorphism with a homoclinic tangency there exist open regions consisting of diffeomorphisms exhibiting nontransversal intersections between stable and unstable manifolds of hyperbolic basic sets. Such sets were called {\em wild hyperbolic} by Newhouse. The original formulation of Newhouse result is as follows:\\

{\bf Newhouse Theorem} \cite{N79}. {\em Let $M^2$ be a $C^\infty$ compact 2-dimensional manifold and let
	$r\geq 2$. Assume that $f \in \mbox{Diff}\, ^r(M^2)$ has a hyperbolic set whose stable and unstable manifolds are tangent at some point $x$. Then $f$ may be $C^r$ perturbed inside an open set
	$U \subset \mbox{Diff}\, ^r(M^2)$ so that each $g \in U$ has a wild hyperbolic set near the orbit of
	$x$.}\\

Several consequences, derived from this theorem, have become crucial in the theory of dynamical systems:
\begin{itemize}
	\item
	There exist open regions in the space of 2-dimensional diffeomorphisms (3-dimensional flows), with the $C^r$-topology, $r\geq 2$, called {\em Newhouse regions}, where the systems having a homoclinic tangency form a dense subset.
	\item
	These Newhouse regions exist in any neighbourhood of any 2-dimensional diffeomorphism having a homoclinic tangency.
\end{itemize}

Newhouse Theorem was extended to a general multidimensional context~\cite{GST93b,PV94,Romero95} and later on to area-preserving diffeomorphisms~\cite{Duarte99,Duarte00,Duarte02}.\footnote{Indeed, it also holds in the multidimensional symplectic case~\cite{DGT}.} 
In the context of general parameter unfoldings~\cite{N79,GST93b}, Newhouse regions are also regarded
as open domains in the parameter space such that the values of the parameters which give rise to homoclinic tangencies form a dense subset. In the case of 1-parameter families, they are usually called {\em Newhouse intervals}.

One of the most known and fundamental dynamical property of Newhouse regions is {\em the coexistence of infinitely many hyperbolic periodic orbits of different types}.  In the dissipative framework, i.e. when the initial quadratic homoclinic tangency is associated to a fixed (periodic) point $O$ with multipliers
$\lambda_1,...,\lambda_m,\gamma$, where $|\lambda_i|<1<|\gamma|$ and the saddle value
$\sigma = \mbox{max}_i |\lambda_i|\cdot |\gamma|$ is less than 1,  this property is known as
{\em Newhouse phenomenon}:
\begin{itemize}
	\item
	In the dissipative case, the set  $\mathcal{ B}$ of parameter values $\mu$ in any Newhouse interval $I$ giving rise to the coexistence of infinitely many periodic sinks and saddles form
	a residual subset.
\end{itemize}

This result was first obtained in~\cite{N74}. Its proof is based essentially on the theory of bifurcations of homoclinic tangencies. The basic elements of this theory were settled in the celebrated work by Gavrilov and Shilnikov \cite{GaS73} where the so-called {\em Theorem on Cascades of periodic sinks} was proved. Indeed, this theorem states the existence of an infinite sequence of intervals of values of a (splitting) parameter for which there exists a single stable periodic.  Multidimensional versions of this result and criteria of birth of periodic sinks at homoclinic bifurcations were established in~\cite{G83,PV94,GStT96,GST08}.

Newhouse phenomenon is very important, in particular, in the theory of the so-called quasiattractors~\cite{AfrSh83}, i.e. strange attractors which either contain periodic sinks of very large periods or such periodic sinks appear under arbitrary small perturbations. Therefore, a natural question arises: how often is the Newhouse phenomenon met in chaotic dynamics?

A partial answer to this question concerning the measure of the set $\mathcal{B}$ introduced above
was considered in a series of papers. Indeed, in~\cite{TY86,Gorodetski-Kaloshin07} the authors showed that this set $\mathcal{B}$ contains a zero-measure secondary subset of parameter values 
for which there exist infinitely many single-round periodic orbits (i.e., orbits passing only once within a neighbourhood of the initial homoclinic orbit).\footnote{
This set $\mathcal{ B}$ can have positive measure for a dense set of suitable families~\cite{T10} and also
for generic families of multidimensional (with $\mbox{dim}\;\geq 3$) diffeomorphisms~\cite{B16}.}

Since Newhouse regions exist near any system presenting a homoclinic tangency, they can be found in the space of parameters of many dynamical models exhibiting chaotic behaviour and in the absence of uniform hyperbolicity. Their extreme richness makes a complete description an unreachable task: tangencies of arbitrarily high order as well as highly degenerate periodic orbits are dense in these regions~\cite{GST93a,GST99}.  
called {\em mixed dynamics} if the closures of the sets of periodic orbits of different types have non-empty intersections. This property can be generic\footnote{It is also persistent in the case of a type of dynamical chaos~\cite{GT17}, which is
characterised by the fundamental property that the intersection of an attractor $\mathcal{A}$ and a repeller $\mathcal{R}$ is non-empty and $\mathcal{ A}\neq \mathcal{ R}$. This is neither the situation in the dissipative chaos (strange attractor), when
$\mathcal{ A}\cap \mathcal{ R}=\emptyset$, nor in the conservative chaos, when $\mathcal{ A} = \mathcal{R}$.}. 
Indeed (see~\cite{GST97}), there exist Newhouse regions with mixed dynamics near any $2$-dimensional diffeomorphism with a nontransversal heteroclinic cycle containing at least two saddle periodic points $O_1,O_2$ whose Jacobians satisfy that $|J(O_1)|>1$ and $|J(O_2)|<1$. This kind of cycles is commonly referred as \emph{contracting-expanding} and it appears to be rather usual in 2-dimensional
reversible diffeomorphisms.

Recall that a diffeomorphism $f$ is called \emph{reversible} if it is smoothly conjugated to its inverse by means of an involution $R$ (named a \emph{reversor}), that is, $R \circ f = f^{-1} \circ R$, with $R^2 = \mathrm{Id}$, $R\ne \mathrm{Id}$. The involution $R$ does not need to be linear. It is often assumed to have the same smoothness as the diffeomorphism $f$. Equivalently, $f$ is reversible if and only if it can be written as the product of two involutions, $f=g \circ h$ with $g^2= h^2 = \mathrm{Id}$. The points which are invariant by the involution $R$ form the symmetry manifold $\mathit{Fix}\,R= \left\{ (x,y) \ | \  R(x,y)=(x,y) \right\}$. Along this work we will consider planar $R$-reversible diffeomorphisms $f$ with $R$ such that $\dim \mathit{Fix}\,R=1$, that is, a curve. 

We say that an object $\Lambda$ is \emph{symmetric} when $R(\Lambda)=\Lambda$. To put more
emphasis, the notation self-symmetric may be used. By a \emph{symmetric couple of objects} $\Lambda_1,\Lambda_2$, we mean two different
objects which are symmetric one to each other, i.e., $R(\Lambda_1)=\Lambda_2$ and, thus, $\Lambda_1=R(\Lambda_2)$.

\begin{figure}[bht]
\begin{center}
\includegraphics[width=12cm]{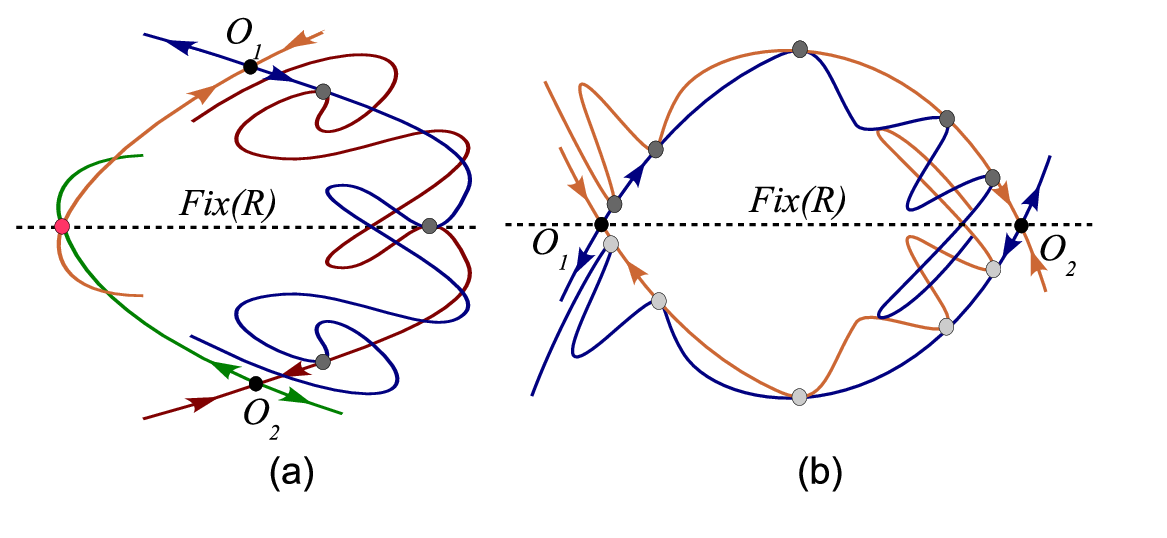}
\caption{{\footnotesize Two different examples of planar reversible maps with symmetric nontransversal
(quadratic tangency) heteroclinic cycles: (a) with a
nontransversal symmetric heteroclinic orbit to a symmetric couple of saddle
points, and (b) with a symmetric couple of nontransversal heteroclinic
orbits to symmetric saddle points.}} \label{fig:Intro1}
\end{center}
\end{figure}

Two examples of contracting-expanding heteroclinic cycle for a $R$-reversible diffeomorphism are shown in Fig.~\ref{fig:Intro1}. In case~(a)
the diffeomorphism has a symmetric couple of saddle periodic (fixed) points~$O_1$ and~$O_2= R(O_1)$, as well as two heteroclinic orbits
$\Gamma_{12}\subset W^u(O_1)\cap W^s(O_2)$ and $\Gamma_{21}\subset W^u(O_2)\cap W^s(O_1)$ such that
$R(\Gamma_{21})=\Gamma_{21}$, $R(\Gamma_{12})=\Gamma_{12}$. The orbit $\Gamma_{12}$ is nontransversal: the manifolds $W^u(O_1)$ and $W^s(O_2)$ have a quadratic tangency along that orbit. Since $R(O_1)=O_2$, their Jacobians verify
$J(O_1)=J^{-1}(O_2)$ and, provided that $J(O_i) \neq \pm 1$, $i=1,2$, the heteroclinic cycle is
contracting-expanding.

Reversible diffeomorphisms can present a very rich dynamics and it is worth studying them by themselves.
Moreover, when they are not conservative (this is an open property) they can possess very interesting dynamics and, in particular, the so-called {\em reversible mixed dynamics}. Its essence, for the $2$-dimensional case, is given by the following two conditions:
\begin{itemize}
\item
The reversible diffeomorphism has simultaneously infinitely many symmetric couples of periodic sinks-sources, periodic saddles with Jacobians greater and less than 1 as well as infinitely many symmetric periodic elliptic orbits and periodic saddles with Jacobian equal to 1.
\item
The closures of periodic orbits of different types have non-empty intersections.
\end{itemize}
These properties seem to be universal when symmetric homoclinic tangencies and symmetric nontransversal heteroclinic cycles are involved in the dynamics. 
Indeed, the following assertion was formulated in~\cite{DGGLS13}.

\medskip

\noindent{\bf Reversible Mixed Dynamics Conjecture (RMD).} \emph{ $2$-dimensional reversible diffeomorphisms with
reversible mixed dynamics are generic in Newhouse regions where diffeomorphisms with symmetric homoclinic or/and
heteroclinic tangencies are dense.}

\bigskip

This RMD conjecture is true when Newhouse regions with $C^r$-topology ($2\leq r \leq \infty$)
are considered (see~\cite{GLRT14}). However, in the real analytic case and for parameter families, it has been proved for a general $1$-parameter unfolding only in two cases -- for 2-dimensional reversible diffeomorphisms with nontransversal heteroclinic cycles, as shown in Fig.~\ref{fig:Intro1}. The cycle of Fig.~\ref{fig:Intro1}(a) was considered in~\cite{LS04}: such
a cycle contains
a symmetric couple of saddle fixed (periodic) points
(with Jacobians less and greater than 1, respectively) and a pair of symmetric transverse and nontransversal heteroclinic orbits. The cycle of Fig.~\ref{fig:Intro1}(b) was considered in~\cite{DGGLS13}: such a cycle contains a symmetric couple of nontransversal heteroclinic orbits to symmetric saddle fixed (periodic) points.

%\iffalse
\begin{figure}[bht]
\begin{center}
\includegraphics[width=12cm]{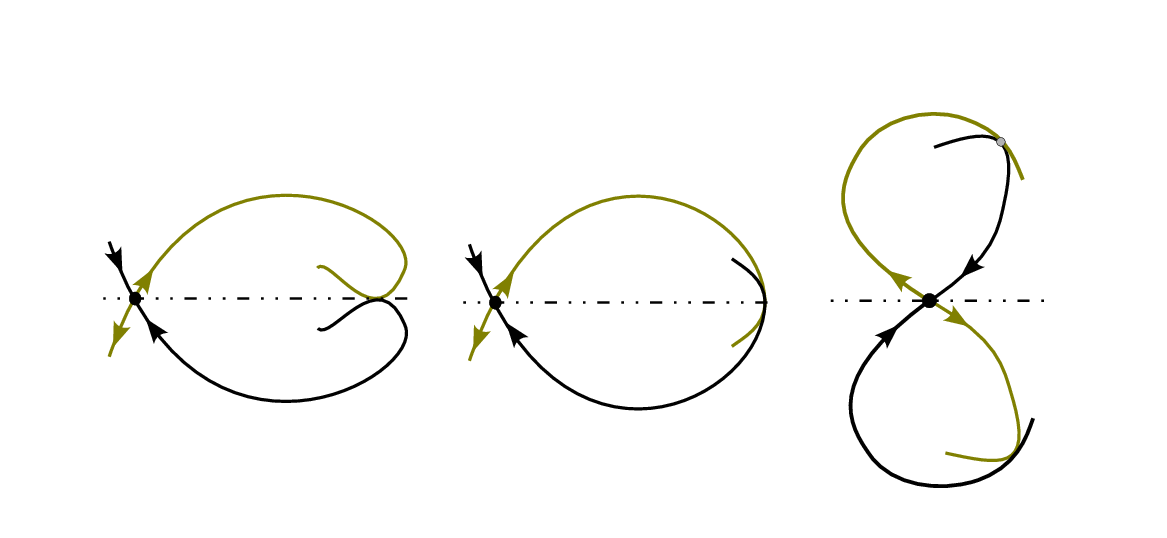}
\caption{{\footnotesize Three examples of planar reversible maps with symmetric nontransversal
homoclinic tangencies: (a) a symmetric quadratic homoclinic tangency;  (b) a symmetric cubic
homoclinic tangency; (c)  a symmetric couple of nontransversal homoclinic (figure-8) orbits to the
same symmetric saddle point.}} \label{fig:3case}
\end{center}
\end{figure}
%\fi

One of the targets concerning RMD conjecture is its proof for 2-dimensional reversible diffeomorphisms which have a homoclinic tangency to a symmetric fixed (periodic) point. There are three main cases, as illustrated in Fig.~\ref{fig:3case}. Figure~\ref{fig:3case}(a) and Figure~\ref{fig:3case}(b) relate to the case when the orbit of the homoclinic tangency is also symmetric and the tangency is either (a)~quadratic or (b)~cubic. In Fig.~\ref{fig:3case}(c) we have the case of a symmetric fixed point and a symmetric couple of orbits with quadratic homoclinic tangencies.

This paper is devoted to this {\em third case} displayed in Fig.~\ref{fig:3case}(c). Roughly speaking, it will be shown that in a general (and symmetrical) unfolding of $1$-parameter families of reversible maps with homoclinic tangencies, there exist Newhouse intervals with reversible mixed dynamics.
%\begin{remark}
We notice that the results of this paper will not only concern orientable planar reversible maps, as the one showed in Fig.~\ref{fig:3case}(c). They will be also valid for maps defined on 2-dimensional non-orientable manifolds allowing a similar structure. For example, on a manifold constructed as a disc surrounding the saddle point with two glued symmetric M\"{o}bius bands.
%\end{remark}

The paper is structured as follows. Section~\ref{sec:mainres} contains the statement of the problem, the main hypotheses and the description of the principal results: Theorems~\ref{th:th1}--\ref{thm:4}. Section~\ref{sec:prelim} deals with the construction of the local and global maps. Theorem~\ref{th:th1} and Theorem~\ref{th:th2} are proved in Section~\ref{sec:proofTh1} and Section~\ref{sec:proofTh2}, respectively.
Section~\ref{sec:proofMainThm} is devoted to the proof of  Theorems~\ref{prop:ni} and~\ref{thm:4}. Finally, 
in Section~\ref{sec:ex} we present some examples of periodically perturbed reversible vector fields giving rise to reversible maps with quadratic hetero/homoclinic tangencies as considered above.

\section{Setting and main results}\label{sec:mainres}

\subsection{The framework}\label{sec:sec_2_1}

Let $f_0$ be a $\mathcal{C}^r$-smooth ($r\geq 4$) reversible
%planar
diffeomorphism of a 2-dimensional manifold $M^2$ with reversor $R$ satisfying $\dim \mathit{Fix}\, R = 1$.
Assume that the following hypotheses hold:
\begin{itemize}
\item[\textsf{[A]}]
The diffeomorphism $f_0$ has a (symmetric) saddle fixed point $O \in \mathit{Fix}\,R$
with multipliers $\lambda,\lambda^{-1}$ and
$0<\lambda<1$.

\item[\textsf{[B]}]
$f_0$ has a symmetric couple of homoclinic orbits $\Gamma_1$ and $\Gamma_2$ such that $\Gamma_2 = R(\Gamma_1)$ (and, thus, $\Gamma_1 = R(\Gamma_2)$)   and satisfies that the invariant manifolds $W^{u}(O)$ and $W^{s}(O)$ have  quadratic tangencies at the points of $\Gamma_{1}$ and $\Gamma_2$.
\end{itemize}
\begin{figure}[bht]
\begin{center}
\includegraphics[width=13cm]{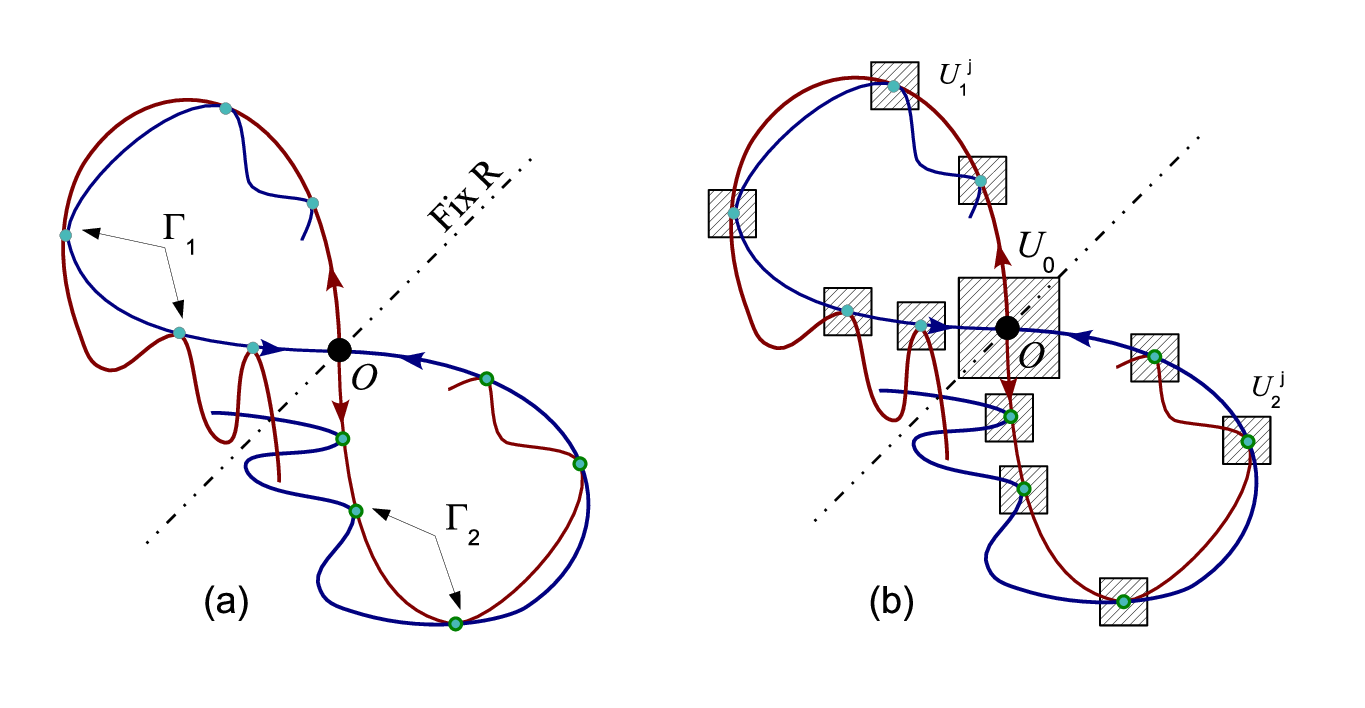}
\caption{(a) An example of reversible map with a couple of symmetric homoclinic tangencies (homoclinic figure-8).
(b) A neighbourhood of the contour $O\cup\Gamma_1\cup\Gamma_2$.}
\label{fig:hom8neigh}
\end{center}
\end{figure}
Let us be more precise with the latter hypothesis.
Take $U$ a small fixed neighbourhood of the contour $O\cup \Gamma_{1}\cup\Gamma_{2}$.
$U$ is formed by the union of a small neighbourhood $U_0$ of the point $O$ and several neighbourhoods $U_1^j$ and $U_2^j$, $j=1,...,n$,  of those points of the orbits $\Gamma_{1}$ and $\Gamma_{2}$ which do not lie in $U_0$ (see Fig.~\ref{fig:hom8neigh}(b)). Thus, $\Gamma_{1}\subset U_0\cup U_1$ and $\Gamma_{2}\subset U_0\cup U_2$, where $U_i = U_0\cup U_i^1...\cup U_i^n$ is a neighbourhood of the homoclinic orbit $\Gamma_i$, for $i=1,2$.
It is not restrictive to assume that $U$ is symmetric, that is $R(U)=U$. Indeed, this comes from assuming that $R(U_0)=U_0$ and $R(U_1^j)=U_2^j$ (and so $R(U_2^j)=U_1^j$).

Consider the orbit $\Gamma_1$ and take a pair of its points, say, $M_1^-\in \Wuloc(O)\cap U_0$ and $M_1^+:=f_0^q(M_1^-) \in \Wsloc(O)\cap U_0$, for a suitable positive integer value $q$.
Denote by $\Pi_1^- \subset U_0$ a small neighbourhood of $M_1^-$ and define the map $T_1:=f_0^q: \Pi_1^- \rightarrow U_0$.
Assume that the following hypothesis is also satisfied:

\begin{itemize}
\item[\textsf{[C]}]
The Jacobian $J_1=J(T_1)|_{M_1^-}$ of the map $T_{1}$ at the point $M_1^-$ is different from $\pm 1$. Without loss of generality we can assume that $|J_{1}|<1$.
\end{itemize}
It is not difficult to check that condition \textsf{[C]} does not depend on the choice of the points $M_1^-$ and $M_1^+$. Moreover, it implies that the map $T_1$, defined in a neighbourhood of $M_1^-$ is not conservative. \\

\begin{figure}[htb]
\begin{center}
\includegraphics[width=14cm]
%, height=7cm]
{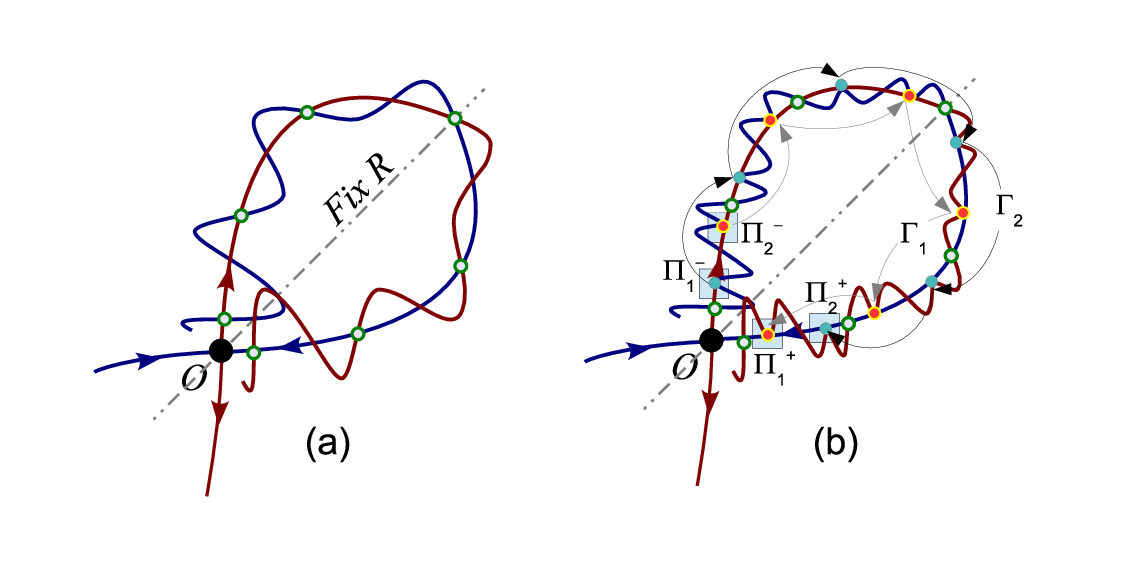} \caption{(a) A reversible diffeomorphism with a symmetric transversal homoclinic
orbit; (b) creation of a symmetric couple of nontransversal homoclinic orbits $\Gamma_1$ and
$\Gamma_2$ (a ``fish'' configuration).}
\label{fig:neigh8rev}
\end{center}
\end{figure}

%\textbf{Remarks.}
\begin{remark}
\mbox{}

\begin{itemize}
	\item[{\bf 1.}] We do not consider the case when the fixed point $O$ has multipliers $\lambda,\lambda^{-1}$ with
		$-1<\lambda<0$. This is a much more complicated case, since $f_0$ would have an additional symmetry due to the negativity of the two multipliers of $O$. 
		%By this reason, for example, the condition C can be fulfilled here.
		\item[{\bf 2.}] In condition \textsf{[C]}, the case $0<J_1<1$ corresponds to $f_0$ orientable while the case $-1<J_1<0$ relates to $f_0$ non-orientable. The latter means that the manifold $M^2$ is non-orientable (the orbit behaves near the global pieces of
		$\Gamma_1$ and $\Gamma_2$, geometrically, like on a M\"{o}bius band).
	%
	%\item
	\item[{\bf 3.}] Our assumptions also cover the case of reversible maps like in Fig.~\ref{fig:neigh8rev}, i.e. when only one pair of stable and unstable separatrices of $O$  create the homoclinic orbits
		$\Gamma_1$ and $\Gamma_2$. Fig.~\ref{fig:neigh8rev}(b) shows how such ``fish'' configuration nontransversal heteroclinic cycle may be created by perturbation of a reversible map with a symmetric transverse homoclinic orbit.
\end{itemize}
\end{remark}

\medskip

Consider two points $M_2^-\in \Wuloc(O)$ and $M_2^+\in \Wsloc(O)$ of the orbit $\Gamma_{2}$ being the symmetric images of the homoclinic points $M_1^+$ and $M_1^-$, i.e.
$M_2^-= R(M_1^+)$ and $M_2^+=R(M_1^-)$. Since $f_0$ is ($R$-)reversible it follows that $f_0^q(M_2^+)=M_2^-$.  Let $T_{2}$ denote the
restriction of the map $f_0^q$ onto a small neighbourhood of the point $M_2^-$. Moreover, we can consider $T_2$ defined from $\Pi_2^- = R(\Pi_1^+)$ onto $\Pi_2^+ = R(\Pi_1^-)$ (see Fig.~\ref{fig:frtms(rev)2}).
Since $T_2 = R(T_1^{-1})$ we have that $J(T_2)|_{M_2^-} = (J(T_1)|_{M_1^-})^{-1}$ and from \textsf{[C]} it follows that $|J_2| = | J(T_2)|_{M_2^-}| >1$.
As it will be properly defined later, iterations of $f_0$ in the neighbourhood $U_0$ around $O$ will be represented by the map $T_0^k$, for positive integer $k$. 

Observe that, for close to $f_0$ maps, one can subdivide nonwandering orbits on $U$ (except for $O$) into three different types: 1-orbits that stay only in $U_0\cap U_1$; 2-orbits  that stay only in $U_0\cap U_2$; and 12-orbits that visit both $U_0\cap U_1$ and $U_0\cap U_2$.
From these types of orbits, we select the so-called single-round periodic orbits, that is those which pass only once inside
$U$. We will refer to them, respectively, as single-round periodic $1$-, $2$- and $12$-orbits.  

For $1$-orbits, we will consider points $x\in \Pi_1^+$, take its image under suitable iterates $k$ of $T_0$, reaching $\Pi_1^-$ and studying 
$\bar x = T_1T_0^k(x)\in \Pi_1^+$ as its return point. If $\bar{x}=x$ we say that $x$ is a fixed point of the first return map $T_{1k} = T_1 T_0^k$.
Analogously, the first return maps for single-round periodic 2-orbits may be represented in the form $T_{2k} = T_2 T_0^k$, from $\Pi_2^+$ onto itself. 
And finally, we will also look for single-round periodic 12-orbits or, equivalently, fixed points of
$\Totkm = T_2 T_0^m T_1 T_0^k$ from $\Pi_2^+$ onto itself, for large enough integers $k$ and $m$. For more details, see Section~\ref{sec:prelim} and Figs.~\ref{fig:hom8neigh} and \ref{fig:frtms(rev)}.

\subsection{The results}\label{sec:sec2_2}

Let $\left\{f_{\mu} \right\}$ be a $1$-parameter family of ($R$-)reversible diffeomorphisms that
unfolds at $\mu=0$ the initial homoclinic tangencies of the diffeomorphism $f_0$ defined above.
Assume that $f_0$ satisfies conditions {\rm \textsf{[A,B,C]}}. Then, the following theorem shows the global symmetry-breaking
bifurcations undergone in this case:
\begin{theorem}
For the family $\left\{f_{\mu} \right\}$, in any segment $[-\mu_0,\mu_0]$ with $\mu_0>0$
small, there are infinitely many intervals $\delta_{k}$,
with boundaries $\mu = \mu_k^+$ and $\mu_k^-$ where $\mu_k^\pm\to 0$ as $k\to\infty$, satisfying:
\begin{itemize}
 \item Symmetric (and simultaneous) single round 1-orbits and 2-orbits of period $k+q$ undergo non-degenerate saddle-node and period-doubling bifurcations at the values $\mu=\mu_{k}^{+}$ and $\mu=\mu_{k}^{-}$, respectively.
 \item The first return maps $T_{1k}$ and $T_{2k}$ have at $\mu\in\delta_{k}$ two fixed points:
a sink and a saddle for $T_{1k}$ and a source and a saddle for $T_{2k}$.
\end{itemize}
\label{th:th1}
\end{theorem}
This theorem can be seen as an extension of the theorem on cascade of periodic sinks in~\cite{GaS73,N74} for the case when the saddle fixed point is conservative and the global dynamics near the homoclinic orbit is dissipative.
In general, these intervals $\delta_{k}$ will be non-intersecting (see Remark~\ref{rmk:non-intersec}
for a wider explanation on that).

In contrast to Theorem~\ref{th:th1}, the following theorem deals with the global bifurcations giving rise to symmetric conservative dynamics, that is, the bifurcations of birth of symmetric single-round elliptic $12$-orbits.
\begin{theorem}
For the family $\left\{f_{\mu} \right\}$ under consideration, in any segment $[-\mu_0,\mu_0]$ with $\mu_0>0$
small, there exist infinitely many intervals $\delta_{km}^c$
accumulating at $\mu=0$ as $k,m\to\infty$ such that the first-return map $\Totkm$ has at $\mu\in\delta_{km}^c$ symmetric elliptic and saddle fixed points.
\label{th:th2}
\end{theorem}
Next result is Newhouse Theorem for the case under consideration. 
%%%Theorem 3
\begin{theorem}
For the family $\left\{f_{\mu} \right\}$, in any segment $[-\mu_0,\mu_0]$ with $\mu_0>0$
small, there exist open intervals $n_i$ such that the set of values $\mu\in n_i$ for which the corresponding map $f_\mu$ satisfying the following two properties \textrm{(a)} and \textrm{(b)} form a dense subset of $n_i$:
\begin{itemize}
 \item[(a)] $f_{\mu}$ has a symmetric couple of homoclinic orbits $\Gamma_{1\mu}\subset U_1$ and $\Gamma_{2\mu}=R(\Gamma_{1\mu})\subset U_2$ to the symmetric saddle fixed point $O_{\mu}$.
 \item[(b)] The manifolds $W^u(O_\mu)$ and $W^s(O_\mu)$ of $f_{\mu}$ have quadratic tangencies at the points of $\Gamma_{1\mu}$ and $\Gamma_{2\mu}$.
\end{itemize}

\label{prop:ni}
\end{theorem}

Summarising, from Theorems~\ref{th:th1}--~\ref{prop:ni} the following result on existence of mixed dynamics is obtained.
%%%%Theorem 4
\begin{theorem}\label{thm:4}
 Let $\{f_\mu\}$ be a $1$-parameter family of $2$-dimensional reversible maps which unfolds at $\mu=0$ a couple of homoclinic tangencies
satisfying conditions \textsf{[A,B,C]}. Then, for any $\mu_0>0$, the
intervals $n_i \subset [-\mu_0,\mu_0]$ from Theorem 3 are Newhouse intervals with reversible mixed dynamics.
\end{theorem}
The proof of Theorems~\ref{th:th1} and~\ref{th:th2} extends along Sections~\ref{sec:prelim}--~\ref{sec:proofTh2}. In contrast, 
the proofs of Theorems~\ref{prop:ni} and~\ref{thm:4} are quite standard and are deferred to the end of the paper:
Theorem~\ref{prop:ni} is proved in Section~\ref{sec:sec6_1} and Theorem~\ref{thm:4} in Section~\ref{sec:sec6_2}.

\section{Preliminary geometric and analytic constructions}
\label{sec:prelim}

\begin{figure}[htb]
\begin{center}
\includegraphics[width=13cm]
%, height=7cm]
{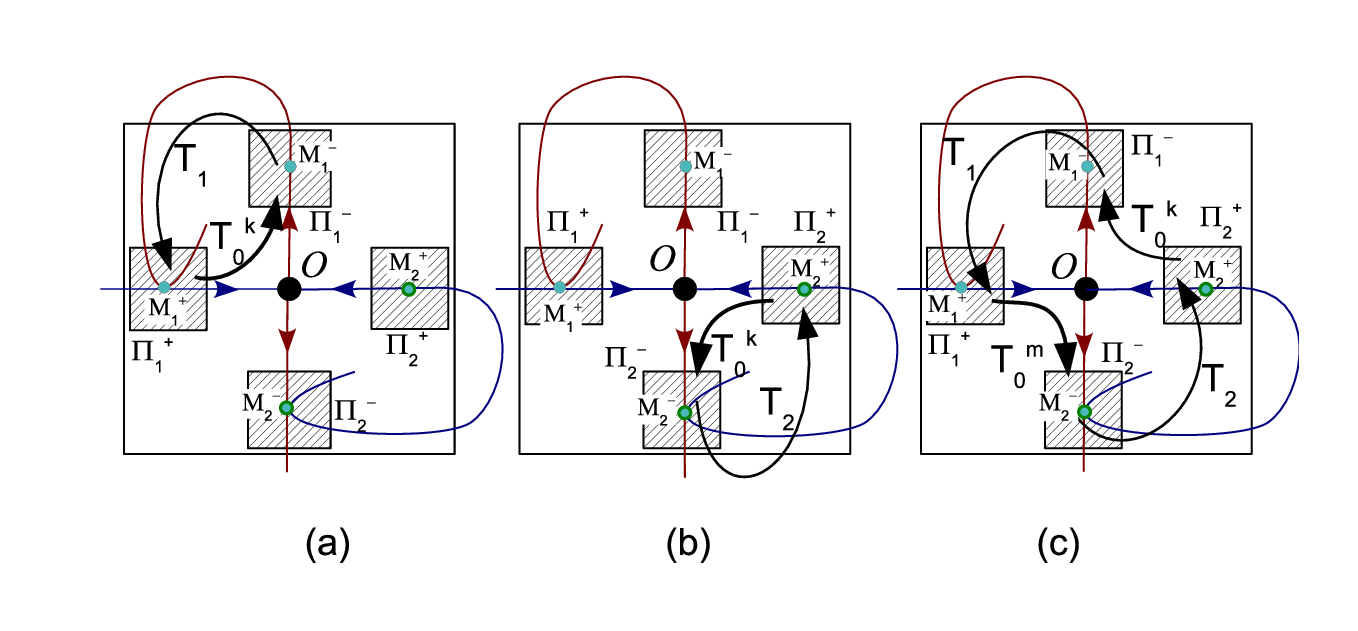} \caption{A geometric structure of the homoclinic points $M_1^+,M_1^-,M_2^+$ and
$M_2^-$ and their neighbourhoods in the case of figure-8 homoclinic configuration. Schematic
actions of the first return maps: (a) $T_{1k}= T_1 T_0^k$, (b) $T_{2k}= T_2
T_0^k$ and (c) $\Totkm= T_2 T_0^m T_1 T_0^k$.  }
\label{fig:frtms(rev)}
\end{center}
\end{figure}

\begin{figure}[htb]
\begin{center}
\includegraphics[width=13.5cm]
%, height=7cm]
{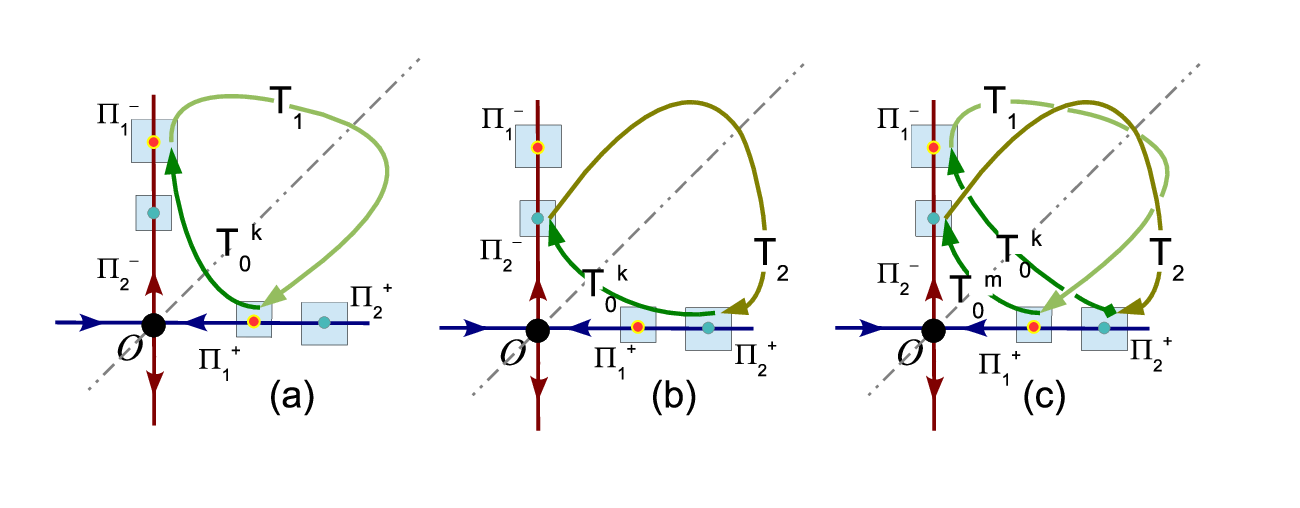} \caption{A geometric structure of the homoclinic points $M_1^+,M_1^-,M_2^+$
and $M_2^-$ and their neighbourhoods in the case of ``fish'' homoclinic configuration. Schematic
actions of the first return maps: (a) $T_{1k}= T_1 T_0^k$, (b) $T_{2k}= T_2
T_0^k$ and (c) $\Totkm= T_2 T_0^m T_1 T_0^k$.}
\label{hom8new1(maps)}
\end{center}
\end{figure}

Let us consider a map $f_\mu$ from our $1$-parameter family and let denote by $T_{0} \equiv {f_\mu}\bigl|_{U_0}$ its restriction onto a neighbourhood $U_0$ of the fixed point $O$. This $\mu$-dependent map $T_{0}$ is called the {\em local map}.
We introduce the so-called {\it global maps} $T_{1}$
and $T_{2}$ through the following relations: $T_{1}\equiv f^q_\mu :\Pi_1^-\to\Pi_1^+$ and $T_{2}\equiv
f^q_\mu :\Pi_2^-\to\Pi_2^+$. They are well defined for small values of $\mu$ since $f^q_0(M_1^-) =
M_1^+$ and $f^q_0(M_2^-) = M_2^+$. Then the {\em first-return maps} $T_{1k}:\Pi_1^+\mapsto\Pi_1^+$,
$T_{2k}:\Pi_2^+\mapsto\Pi_2^+$ and $\Totkm:\Pi_1^+\mapsto\Pi_1^+$ can be defined by the following
composition of maps:
\begin{equation*}
\begin{array}{l}
\Pi_1^+ \;\stackrel{T_{0}^k}{\longrightarrow}\;\;\Pi_1^- \;\stackrel{T_{1}}{\longrightarrow}\;\; \Pi_1^+\;, \\

\Pi_2^+ \;\stackrel{T_{0}^k}{\longrightarrow}\;\;\Pi_2^- \;\stackrel{T_{2}}{\longrightarrow}\;\; \Pi_2^+ \;, \\

\Pi_2^+ \;\stackrel{T_{0}^k}{\longrightarrow}\;\;\Pi_1^- \;\stackrel{T_{1}}{\longrightarrow}\;\;
\Pi_1^+ \;\stackrel{T_{0}^m}{\longrightarrow}\;\;\Pi_2^-\;\stackrel{T_{2}}{\longrightarrow}\;\;
\Pi_2^+\;,
\end{array}
%\label{TkmPi}
\end{equation*}
(see Fig.~\ref{fig:frtms(rev)} and~\ref{hom8new1(maps)}). In short, we will denote these compositions by $T_{1k}= T_1 T_0^k$, $T_{2k}=T_2 T_0^k$ and $\Totkm= T_2 T_0^m T_1 T_0^k$.
As it is usual in this kind of problems, one seeks for suitable local coordinates on $U_0$ in which the map $T_{0}$ exhibits its simplest form. The following lemma introduces
$C^{r-1}$-coordinates that allow our local map $T_0$ to be written in the so-called {\em (saddle) normal
form} or first order (saddle) normal form.
\begin{lm}[Saddle Normal Form~{\cite[Lemma 1]{GS90}}]
Assume $r\geq 4$
and let $T_0$ be a $C^r$-smooth reversible planar map with reversing (nonlinear in general) involution $R$ satisfying that $\dim\mathrm{Fix}\,R = 1$.
Suppose that $T_0$ has a saddle fixed (periodic) point $O$ at the origin which belongs to
$\mathrm{Fix}\,R$ and has multipliers $\lambda$ and $\lambda^{-1}$, with $|\lambda|<1$. Then
there exist $C^{r-1}$-smooth local coordinates near $O$ in which the map $T_0$
%(or $T_0^n$, where $n$ is the period of $O$) This could be written, if necessary, as a posterior remark
can be written in the so-called Shilnikov cross-form:
\begin{equation}
T_0:\;\; \left\{
\begin{array}{rcl}
\bar x &=& \lambda x  + h(x,\bar y)x^2\bar y,\\
y &=& \lambda \bar y + h(\bar y,x)x\bar y^2.
\label{CrossFormSaddle}
\end{array}
\right.
\end{equation}
\label{LmSaddle}
\end{lm}

\begin{remark}
In these local coordinates the map $T_0$ is reversible under the linear involution $L(x,y)=(y,x)$.
Indeed (see~\cite{DGGLS13}, for instance), it is enough to check that $(L  T_0  L )^{-1} = T_0$. Observe that
\[
L T_0  L:
\left\{
\begin{array}{rcl}
\bar y &=& \lambda y  + h(y,\bar x)y^2 \bar x,\\
x &=& \lambda \bar x + h(\bar x,y)y\bar x^2
\end{array}
\right.
\]
and thus $(L  T_0  L )^{-1}$,  which corresponds to interchange $(\bar x, \bar y)\leftrightarrow (x,y)$, gives rise to the expression for $T_0$. Bochner theorem~\cite{MZ55} ensures the simultaneous conjugation of both the map and the reversor.
\end{remark}

\noindent Next lemma provides a suitable expression for the iterates of $T_0$. Namely,
\begin{lm}[{see~\cite{GS90}}] %[Iterations of the local map]
Let $T_0$ be a $C^r$-smooth $R$-reversible map written in (local) normal form~(\ref{CrossFormSaddle}) in a neighbourhood $V$ of a saddle fixed point $O$. Let us consider iterates of $T_0$ in $V$:  $(x_0,y_0),\dots,(x_{\ell},y_{\ell})$ such that $(x_{\ell+1},y_{\ell+1})= T_0(x_{\ell},y_{\ell})$, $\ell=0,\dots,j-1$. Then, one has that
\begin{eqnarray}
x_j &=& \lambda^j x_0 \left(1 + j \lambda^j h_j(x_0,y_j)\right), \label{LocalMapk} \\
y_0 &=& \lambda^j y_j \left(1 + j
\lambda^j h_j(y_j, x_0)\right), \nonumber
\end{eqnarray}
where
the functions  $h_j(x_0,y_j)$ are $\mathcal{O}_2(x_0,y_j)$
and satisfy that they and all their derivatives up to order $r-2$ are
uniformly bounded with respect to $j$ . \label{LmLocalMap}
\end{lm}

Lemmas~\ref{LmSaddle} and~\ref{LmLocalMap} are also valid if $T_0$ depends on parameters. Moreover, if $T_0$ is $C^r$ with respect to both coordinates and parameters,
it can be seen that the normal form (\ref{CrossFormSaddle}) is $C^{r-1}$ with
respect to the coordinates and $C^{r-2}$ with respect to the parameters. Moreover, the derivatives in~(\ref{LocalMapk}) with respect to the parameters and up to order $r-2$
have order $O\left((\lambda+\epsilon)^j\right)$ for any $\epsilon>0$ (we refer the reader to~\cite{GST08}, Lemmas 6 and 7, for more details).

\subsection{Construction of the local and global maps}

We choose in $U_0$ the local coordinates  $(x,y)$ given in Lemma~\ref{LmSaddle}. In these coordinates,
the local stable and unstable invariant manifolds of the point $O$ are straightened: 
$\Wuloc(O)$ can be represented by $x=0$ and $\Wsloc(O)$ by $y=0$.
Moreover, the previously chosen homoclinic points read as follows: $M_1^+=(x_1^+,0)$,
$M_1^-=(0,y_1^-)$, $M_2^+=(x_2^+,0)$ and $M_2^-=(0,y_2^-)$. Since $R(M_1^+) = M_2^-$ and
$R(M_1^-)=M_2^+$, we have that they are $L$-symmetric and therefore $x_1^+=y_2^-$ and $y_1^-=x_2^+$. From the geometry of the figure-8 homoclinic case (see Fig.~\ref{fig:frtms(rev)}) we can assume that
\begin{equation}
x_1^+= y_2^- = - \alpha<0,\qquad y_1^- = x_2^+ = \beta >0
\label{xeqy}
\end{equation}
Analogously, in the ``fish'' configuration we have that
%$x_1^+= y_2^- >0,y_1^- = x_2^+ >0$
$\alpha <0$ and $\beta>0$ (see Fig.~\ref{hom8new1(maps)}).

\begin{figure}[htb]
\begin{center}
\includegraphics[width=14cm]%, height=7cm]
{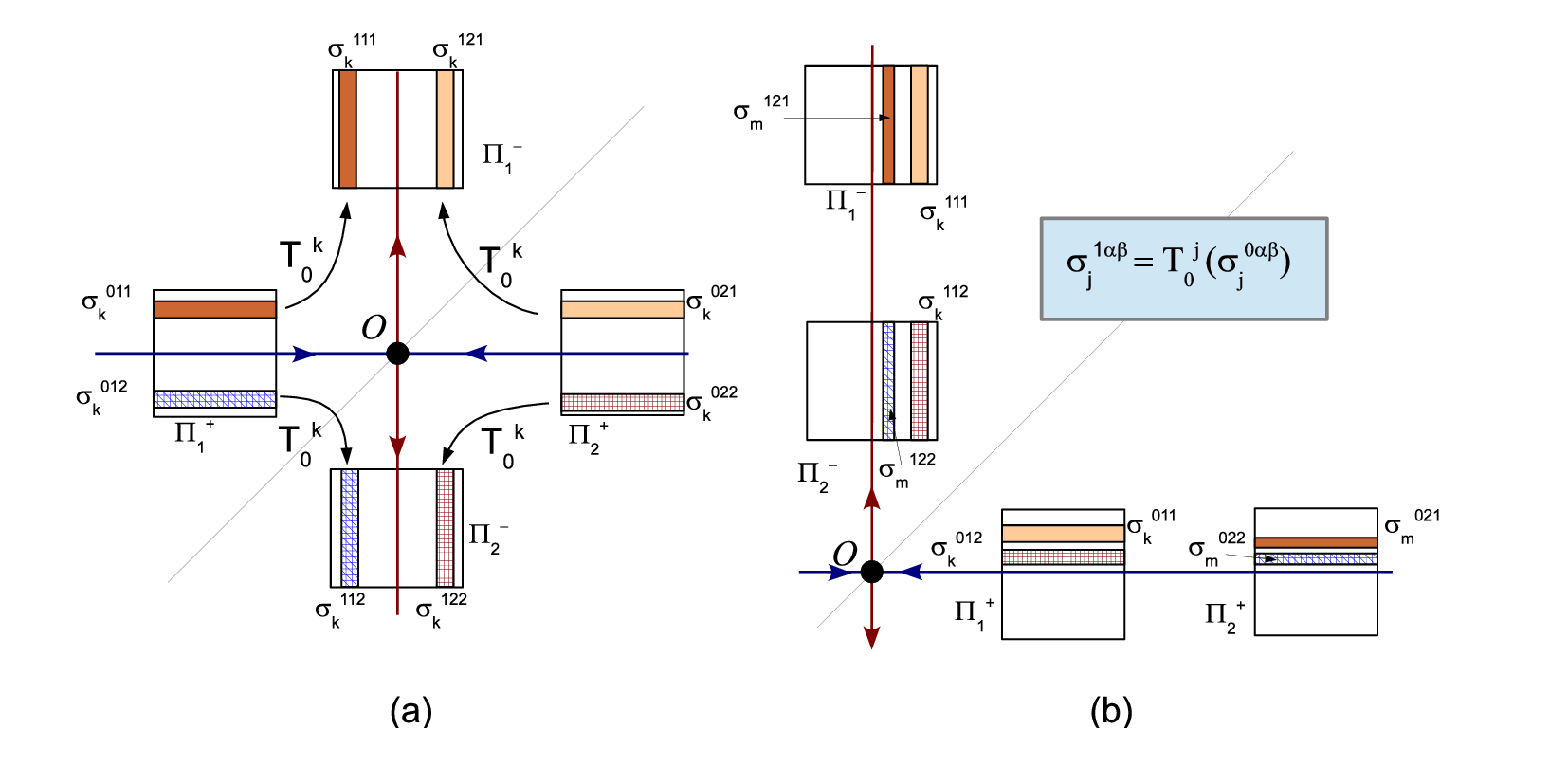} \caption{Domains of definition and range  of the successor map from $\Pi_i^{+}$
into $\Pi_j^-$, $i,j = 1,2$, under iterations of $T_{0}$ in the cases of (a) homoclinic figure-8;
(b) homoclinic ``fish''.  }
\label{fig:locmaps}
\end{center}
\end{figure}

\begin{figure}[htb]
\begin{center}
\includegraphics[width=12cm]
%, height=7cm]
{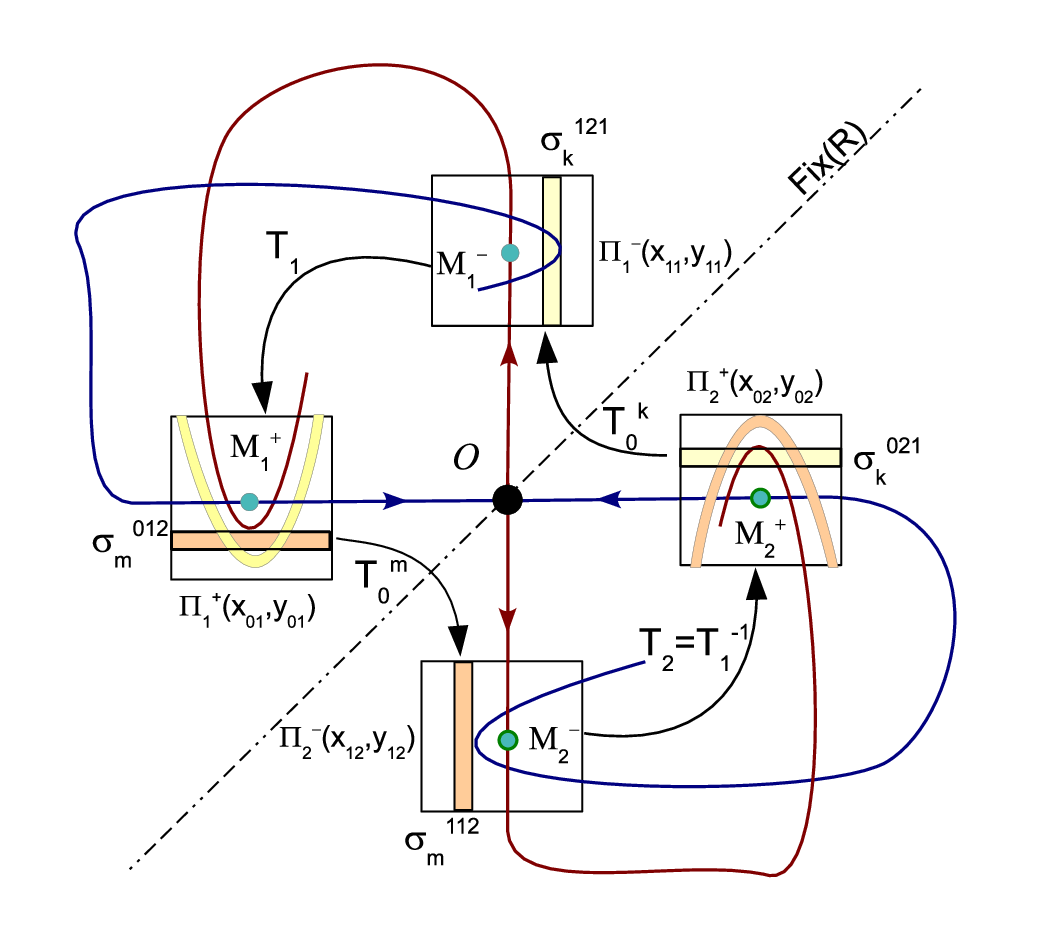} \caption{Domains of definitions and associated coordinates for the first return map $\Totkm = T_2T_0^mT_1T_0^k$.}
\label{fig:frtms(rev)2}
\end{center}
\end{figure}

It is not restrictive to assume that $T_{0}(\Pi_i^-)\cap\Pi_i^-=\emptyset$,
$i=1,2$ (if not, one can reduce the size of $\Pi_i^-$). Therefore, the domains of definition of the
transfer map from $\Pi_i^{+}$ into $\Pi_j^-$, $i,j = 1,2$,  under iterations of $T_{0}$ consist
of infinitely many non-intersecting strips $\sigma_{k}^{0ij}$ which belong to $\Pi_i^+$ and
accumulate at $\Wsloc(O)\cap\Pi_i^+$ as $k\to\infty$.
On its turn, the range of the transfer map consists of infinitely many  strips $\sigma_{k}^{1ij}=
T_{0}^k(\sigma_k^{0ij})$ belonging to $\Pi_i^-$ and accumulating at $\Wuloc(O)\cap\Pi_i^-$ as
$k\to\infty$ (see Figure~\ref{fig:locmaps}).
So, our first return maps are defined on those strips in the following way:
\begin{equation*}
\begin{array}{l}
T_{1k} = T_1 T_0^k : \sigma_k^{011}\stackrel{T_{0}^k}{\longmapsto}
\sigma_k^{111}\stackrel{T_1}{\longmapsto} \sigma_k^{011} \;, \\
T_{2k} = T_2 T_0^k : \sigma_k^{022}\stackrel{T_{0}^k}{\longmapsto}
\sigma_k^{122}\stackrel{T_2}{\longmapsto} \sigma_k^{022} \;, \\
\Totkm = T_2 T_0^m T_1 T_0^k : \sigma_k^{021}\stackrel{T_{0}^k}{\longmapsto}
\sigma_k^{121}\stackrel{T_1}{\longmapsto} \sigma_m^{012}\stackrel{T_{0}^m}{\longmapsto}
\sigma_m^{112}\stackrel{T_2}{\longmapsto} \sigma_k^{021} \;.
\end{array}
%\label{Tkmxy}
\end{equation*}
For large enough values of $k$, Lemma~\ref{LmLocalMap} asserts that the map $T_{0}^k:\sigma_k^{0ij}\left\{(x_0,y_0)\right\}\mapsto\sigma_k^{1ij}\left\{(x_k,y_k)\right\}$ can be written in the form
\begin{equation}
T_{0}^k\;:\; \left\{
\begin{array}{l}
x_{k} = \lambda^k x_{0} \left(1 + k\lambda^k h_k(x_{0},y_{k}) \right),\\
y_{0} = \lambda^k y_{k} \left(1 + k\lambda^k h_k(y_{k},x_{0}) \right)
\end{array}
\right. \label{T1k}
\end{equation}
where $(x_0,y_0)\in \sigma_k^{0ij},\; (x_1,y_1)\in \sigma_k^{1ij}$, $i,j=1,2$. In the ``fish'' configuration case this corresponds to $T_0^k:\Pi_1^+ \left\{(x_{01},y_{01})\right\} \mapsto
\Pi_1^- \left\{(x_{11},y_{11})\right\}$ while in the figure-8 situation this becomes $T_0^k:\Pi_2^+ \left\{(x_{02},y_{02})\right\} \mapsto \Pi_1^- \left\{(x_{11},y_{11})\right\}$ and
$T_0^m:\Pi_1^+ \left\{(x_{01},y_{01})\right\} \mapsto \Pi_2^- \left\{(x_{12},y_{12})\right\}$ (see Fig.~\ref{fig:frtms(rev)2}).
The global map $T_{1}:\Pi_1^-\to\Pi_1^+$  admits the following form
\begin{equation}
\fl T_{1}: \left\{
\begin{array}{rl}
x_{01} - x_1^+ &= F_{1}(x_{11},y_{11}-y_1^-,\mu) \\
               & \equiv a x_{11} + b (y_{11}-y_1^-) + \varphi_1(x_{11},y_{11},\mu),\\
y_{01} &= G_{1}(x_{11},y_{11}-y_1^-,\mu) \\
       &\equiv  \mu + c x_{11} + d (y_{11}-y_1^-)^2  + \varphi_2(x_{11},y_{11},\mu),
\end{array}
\right.
\label{T12}
\end{equation}
where $F_{1}(0)=G_{1}(0)=0$  (since $T_{1}(M_1^-) = M_1^+$ at $\mu=0$) and
$\varphi_1=\mathcal{O}\left((y_{11}-y_1^-)^2 + x^2_{11}\right)$,
$\varphi_2=\mathcal{O}\left( x_{11}^2 + |y_{11}-y_1^-|^3 + |x_{11}||y_{11}-y_1^-|\right).$
%\label{T12hot}
Since $\Wuloc(O))$ and $\Wsloc(O)$ have (local) expressions $\{x_{11}=0\}$ and $\{y_{01}=0\}$
and $T_{1}(\Wuloc(O))$ and $\Wsloc(O)$ undergo a quadratic tangency at $\mu=0$, this implies that
\[
\frac{\partial G_{1}(0)}{\partial y_{11}} = 0,\;\; \frac{\partial^2 G_{1}(0)}{\partial y_{11}^2} =
2d \neq 0.
\]
Its Jacobian $J(T_{1})$ has the form
\begin{equation}
J(T_{1})= -bc + \mathcal{O}\left(|x_{11}| + |y_{11}-y_1^-|\right),
\label{detT12*}
\end{equation}
and so $J_1=J(T_1)|_{M^-} = - bc$ where $0<|bc|<1$ by condition \textsf{[C]}.

Concerning the global map $T_{2}$, its expression is closely related to that of $T_1$. Indeed,
reversibility implies that
$T_{2} = R\; T_{1}^{-1}\;R$ or, equivalently, $T_{1} = R\; T_{2}^{-1}\;R$.
Then, by expression~(\ref{T12}) and having in mind the local $L$-reversibility on the domains $\Pi^-_2$ (Bochner's theorem ensures its conjugation with the non-linear reversor $R$)
we obtain that the map
$T_{2}^{-1}:\Pi_2^+\{(x_{02},y_{02})\}\mapsto \Pi_2^-\{(x_{12},y_{12})\}$ can be written as
\begin{equation*}
\fl T_{2}^{-1}: \left\{
\begin{array}{rl}
x_{12} =& G_{1}(y_{02},x_{02}-x_2^+,\mu) = \\
       &\mu + c y_{02} + d (x_{02}-x_2^+)^2 +  \varphi_2(y_{02},x_{02},\mu),\\
y_{12} - y_2^- =& F_{1}(y_{02},x_{02}-x_2^+,\mu) =  \\
 &a y_{02} + b (x_{02}-x_2^+) + \varphi_1(y_{02},x_{02},\mu),
\end{array}
\right.
%\label{T21-}
\end{equation*}
which means to write $x_1^+=y_2^-$, $y_1^-=x_2^+$ in~(\ref{T12}) and to swap $x\leftrightarrow y$ variables, i.e. $x_{01} \leftrightarrow y_{12}$ and $x_{11} \leftrightarrow y_{02}$.
%$\square_{01} \leftrightarrow \square_{12}$ and $\square_{11} \leftrightarrow \square_{02}$.
As it was done in a previous remark, this expression defines the map
$T_{2}:\Pi_2^-\{(x_{12},y_{12})\}\mapsto\Pi_2^+\{(x_{02},y_{02})\}$ in the implicit form:
$x_{12} = G_{1}(\bar y_{02},\bar x_{02}-x_2^+,\mu), y_{12} - y_2^- = F_{1}(\bar y_{02},\bar
x_{02}-x_2^+,\mu)$ by swapping bar and no-bar variables.

\section{Proof of Theorem~\ref{th:th1}}
\label{sec:proofTh1}
This proof is mainly based on Lemma~\ref{lm:rescMar} which provides, by computing the corresponding equations
and performing a suitable rescaling, an asymptotic expression for the first return map
for large enough values of $k$. Rescaling method has become, since the work of  Tedeschini-Yorke~\cite{TY86},
a very useful tool when dealing with homoclinic connections (see also~\cite{GStT96,GG04,GGT07,GST08,GGS10} and references therein for many examples of such use).
\begin{lm}
\label{lm:rescMar}
Let $\left\{f_\mu \right\}$ be the family under consideration satisfying conditions~\textrm{\textsf{[A,B,C]}}.
Then, for large enough values of $k$, the first return map $T_{1k} :
\sigma_k^0 \rightarrow \sigma_k^0$ can be brought, by a linear change of coordinates and
a convenient rescaling, to the following form
\begin{equation}
\label{eq:FRMmar}
\bar{X}=Y + k\lambda^k h_k^1(X,Y), \qquad \bar{Y} = M_1 + M_2 X - Y^2 + k\lambda^k h_k^2(X,Y), %\mathcal{O}(\lambda^k),
\end{equation}
with
\begin{equation}
\label{eq:M1M2}
%\[
M_1= -d\lambda^{-2k}\left(\mu  - \lambda^k(y_1^- - cx^+) + \tilde\rho_k\right),
%[ \alpha + \beta)+ \mathcal{O}(\lambda^{-k}),
\qquad M_2=bc,
\end{equation} %\]
where $\tilde\rho_k = \mathcal{O}(k\lambda^{k})$ is a small constant and the functions $h_k^j$ have all
their derivatives uniformly bounded up to order $(r-2)$.
\end{lm}

{\em Proof}.
To ease its reading we give first a ``lightweight'' proof of the lemma for a simpler case, i.e. when
the local map $T_0$ is linear, $\bar x = \lambda x, \bar y = \lambda^{-1} y$, and the global map has the form:
\[
\begin{array}{l}
\bar x_0 - x^+ = a x_1 + b(y_1-y^-),\\
\bar y_0 = \mu + c x_1 + d(y_1-y^-)^2 + f_{11} x_1(y_1-y^-).
\end{array}
\]
We have only considered linear terms in the first equation and up to quadratic terms in the second one. We use also (only for a simplification of formulas) the notation
$x^+=x_1^+, y^- = y^-_1$ and denote the coordinates on $\Pi_1^+$ as $(x_0,y_0)$  and on $\Pi_1^-$ as
$(x_1,y_1)$. Then the first return map $T_{1k} = T_1T_0^k$ is written as
\[
\begin{array}{l}
\bar x_0 - x^+ = a \lambda^k x_0 + b(y_1-y^-),\\
\lambda^k\bar y_1 = \mu + c \lambda^k x_0 + d(y_1-y^-)^2 + f_{11} \lambda^k x_0(y_1-y^-),
\end{array}
\]
This (first) highly simplified case will serve the reader (we hope) to be familiar with the different transformations we apply to get the asymptotic H\'enon  map. The general case (that is included rear after this one) will follow the same ideas and procedure.

Introduce the coordinates $\xi = x_0 - x^+, \eta = y_1 - y^-$. Then $T_{1k}$ reads
\begin{equation}
\begin{array}{l}
\bar\xi = a \lambda^k \xi + b\eta + a \lambda^k x^+ ,\\
\lambda^k\bar\eta  = m_1 + c \lambda^k \xi + d\eta^2 + f_{11} \lambda^k \xi\eta +
f_{11} \lambda^k x^+ \eta,
\end{array}
\label{T1xieta}
\end{equation}
%$$
where $m_1 = \mu + c\lambda^k x^+ - \lambda^k y^-$. 

Further, we make one more coordinate shift, $\xi = x + \alpha_k, \eta = y + \beta_k$ with small coefficients $\alpha_k =\mathcal{O}(\lambda^k)$ and $\beta_k = \mathcal{O}(\lambda^k)$, in order to vanish the constant terms in the first equation and the linear in $y$ terms in the second one. Then we obtain
\[
\begin{array}{l}
\bar x = a \lambda^k x + b y + \left[b\beta_k - \alpha_k  + a\lambda^k x^+ + a\lambda^k\alpha_k \right] ,\\
\lambda^k\bar y  = m_2   + (c + f_{11}\beta_k) \lambda^k x  + d y^2 + f_{11} \lambda^k xy +
 \left(2d\beta_k + f_{11} \lambda^k x^+  + f_{11}\lambda^k\alpha_k \right)  y,
\end{array}
\]
where $m_2 = m_1 + \lambda^k(c\alpha_k - \beta_k + f_{11}x^+\beta_k + f_{11}\alpha_k\beta_k) + d\beta_k^2 =
m_1 + \mathcal{O}(\lambda^{2k})$.
The expressions in square brackets are nullified at
\begin{equation}%$$
\alpha_k = \left(ax_1^+ - \frac{bf_{11}x^+}{2d}\right) \lambda^k + \mathcal{O}(\lambda^{2k}), \;\;
\beta_k =  - \frac{f_{11}x^+}{2d} \lambda^k + \mathcal{O}(\lambda^{2k}).
\label{ab}
\end{equation}%$$
For such choice of $\alpha_k$ and $\beta_k$, the map $T_{1k}$ takes the form
\[
\begin{array}{l}
\bar x = a \lambda^k x + b y ,\\
\bar y  = \lambda^{-k}m_2   + (c + \phi_k) x  + d \lambda^{-k} y^2 + f_{11} xy,
\end{array}
\]
where $\phi_k = \mathcal{O}(\lambda^k)$ is a small coefficient. Now, by rescaling the coordinates,
\begin{equation} %$$
x = - \frac{b}{d} \lambda^k X,\;\; y = - \frac{1}{d} \lambda^k Y,
\label{rescT1k}
\end{equation} %$$
we bring the map $T_{1k}$ to the claimed form:
\[
\bar X = Y + \mathcal{O}(\lambda^{k}), \;\; \bar Y = M + bc X - Y^2 + \mathcal{O}(\lambda^{k}),
\]
where $M = -d\lambda^{-2k}m_2 = -d\lambda^{-2k} \left[\mu + (c  x^+_1 - y^-_1) \lambda^k + \mathcal{O}(\lambda^{2k}) \right]$.

\vskip1cm

Let us now deal with the \emph{general case}, that is, with $T_0^k$ given by
\begin{eqnarray*}
x_k &=& \lambda^k x_0 \left( 1 + k\lambda^k h_k(x_0,y_k) \right) \\
y_0 &=& \lambda^k y_k \left( 1 + k\lambda^k h_k(y_k,x_0) \right)
\end{eqnarray*}
and the global map $T_1$ given by
\begin{eqnarray*}
\bar{x}_0 - x^+ &=& a x_1 + b(y_1 - y^-) + \mathcal{O}\left( (y_1-y^-)^2, x_1^2, (y_1-y^-)x_1 \right), \\
\bar{y}_0 &=& \mu + c x_1 + d(y_1-y^-)^2 + f_{11} x_1 (y_1-y^-) + \mathcal{O}\left( x_1^2,(y_1-y^-)^3 \right).
\end{eqnarray*}
Consider the map $T_{1k}=T_1 T_0^k$ and apply the change of coordinates: $\xi=x_0-x^+$, $\eta= y_k - y^-$. Then, $T_{1k}$ admits the expression
\begin{eqnarray}
\bar{\xi} &=& a \lambda^k \xi + b\eta + \left( \lambda^k a x^+ + \mathcal{O}(k \lambda^k) \right) + \gamma_1 \eta^2 + \gamma_2 \lambda^k \xi \eta + \lambda^k \eta, \label{eq:gencase:xieta}\\
\lambda^k \bar{\eta} (1+\mathcal{O}(k\lambda^k))&=& \left( \mu + c\lambda^k x^+ + c\lambda^k(\xi + x^+) k \lambda^k h_k + f_{11} k \lambda^{2k} (\xi + x^+) \eta h_k + \right. \nonumber\\
&& \left. \gamma_1 \lambda^{2k} (\xi+x^+)^2 (1+k\lambda^k h_k) %(^2 + \gamma_2 \eta^3) 
\right) + c\lambda^k \xi + d\eta^2 +   \\
&& f_{11} \lambda^k \xi \eta + f_{11} \lambda^k x^+ \eta. \nonumber
\end{eqnarray}
Following the same steps as for the simplified case, we consider the following \emph{shift}:
\[
\xi=x + \alpha_k, \qquad 
\eta= y + \beta_k 
\]
with $\alpha_k,\beta_k$ to be determined in such a way that the constant term in the equation for $\bar{x}$ and the coefficient of $y$ in $\bar{y}$ both vanish. After performing this shift, equations~\eqref{eq:gencase:xieta} become
\begin{equation}
\label{eq:gencase:xbar}
\bar{x}=a \lambda^k x + by + \left( (a\lambda^k -1) \alpha_k + b(1+\lambda^k) \beta_k + \lambda^k a x^+ + \mathcal{O}(k\lambda^k) \right)
\end{equation}
and
\begin{eqnarray}
\label{eq:gencase:ybar}
\lefteqn{\lambda^k \bar{y} = \left( \mu + c \lambda^k (x^+ - y^-) + ck \lambda^{2k} (\alpha_k + x^+) h_k^0 + f_{11} k \lambda^{2k} (\alpha_k + x^+) \beta_k h_k^0 +  \right. } \nonumber \\
&& \left. \gamma_1 \lambda^{2k} (\alpha_k + x^+)^2
+ \gamma_2 \beta_k^3 + c \lambda^k \alpha_k + d \beta_k^2 + f_{11} \lambda^k \alpha_k \beta_k + f_{11} \lambda^k x^+ \beta_k - \lambda^k \beta_k +   \right. \nonumber \\
&& 
\left. \mathcal{O}(k \lambda^{4k} )\right)
+ \left( c \lambda^k + f_{11} \lambda^k \beta_k + c k \lambda^{2k} h_k^0 + f_{11} k \lambda^{2k} \beta_k h_k^0 \right) x +  \\
&& \left( f_{11} \lambda^k (1+k\lambda^k h_k^0) \alpha_k + 2d \beta_k + 3 \gamma_2 \beta_k^2 + f_{11} k \lambda^{2k} x^+ h_k^0 + f_{11} \lambda^k x^+ \right) y +  \nonumber \\
&& 
\left( d+ 3 \gamma_2 \beta_k \right) y^2 +
\left( f_{11} k \lambda^{2k} h_k^0 + f_{11} \lambda^k \right) xy + \mathcal{O}_3(x,y), \nonumber
\end{eqnarray}
where $h_k^0$ stands for the constant term of $h_k(x^+ + \xi,y^- + \eta)$ 
in $(\xi,\eta)$-variables and we have taken into account that $(1+\mathcal{O}(k\lambda^k))^{-1}= \mathcal{O}(k\lambda^k)$. Thus, we determine $\alpha_k,\beta_k$ to satisfy
\begin{eqnarray}
\label{eq:gencase:alphak:betak}
(a\lambda^k -1) \alpha_k + b(1+\lambda^k) \beta_k &=& - \lambda^k a x^+ + \mathcal{O}(k\lambda^k) \\
f_{11} \lambda^k (1+k\lambda^k h_k^0) \alpha_k + 2d \beta_k + 3 \gamma_2 \beta_k^2 &=& - f_{11} k \lambda^{2k} x^- h_k^0 - f_{11} \lambda^k x^+.
\nonumber
\end{eqnarray}
It is straightforward to check that $\alpha_k, \beta_k =\mathcal{O}(\lambda^k)$. Now, consider the linear system
\begin{eqnarray*}
%\label{eq:gencase:IFTsys}
(a\lambda^k -1) \alpha_k + b(1+\lambda^k) \beta_k &=& - \lambda^k a x^+ + \mathcal{O}(k\lambda^k) \\
f_{11} \lambda^k (1+k\lambda^k h_k^0) \alpha_k + 2d \beta_k  &=& - f_{11} k \lambda^{2k} x^- h_k^0 - f_{11} \lambda^k x^+.
\end{eqnarray*}
This linear system has solutions 
\begin{equation}
\label{eq:gencase:alphabetak}
\begin{array}{rcl}
\alpha_k^0 &=& {\ds \left( a x^+ + \frac{b f_{11}}{2 d} \right) \lambda^k + \mathcal{O}(k \lambda^k) } \\
\beta_k^0  &=& {\ds - \frac{f_{11} x^+}{2d} \lambda^k + \mathcal{O}(k \lambda^k)}.
\end{array}
\end{equation}
Since $d\neq 0$, the determinant
\[
\left| \begin{array}{cc}
             a \lambda^k -1 & b (1+\lambda^k) \\
             f_{11} \lambda^k + \mathcal{O}(k\lambda^{2k}) & 2 d 
            \end{array}
\right| =  - 2d + (2 ad - b f_{11}) \lambda^k - b f_{11} \lambda^{2k} + \mathcal{O}(k \lambda^{2k}) \neq 0,
\]
and so by the Implicit Function Theorem, there exist $\alpha_k=\alpha_k^0 + \mathcal{O}(k\lambda^k)$ and $\beta_k=\beta_k^0 + \mathcal{O}(k\lambda^k)$ solutions of~\eqref{eq:gencase:alphak:betak}, which are $\mathcal{O}(k\lambda^k)$-close to $\alpha_k^0,\beta_k^0$. Thus, considering the shift $\xi=x+\alpha_k$, $\eta=y + \beta_k$, with these already determined $\alpha_k,\beta_k=\mathcal{O}(\lambda^k)$, one gets the following equations for $T_{1k}$:
\begin{equation}
\begin{array}{rcl}
\bar{x} &=& a \lambda^k x + by + \gamma_1 y^2 \\
\lambda^k \bar{y} &=& m_2 + (c\lambda^k + \mathcal{O}(\lambda^{2k}) x + (d+\mathcal{O}(\lambda^k)) y + (f_{11}\lambda^k + \mathcal{O}(k \lambda^{2k})) xy + \mathcal{O}_3(x,y), 
\end{array}
\end{equation}
where
$m_2:= \mu + c \lambda^k x^+ - \lambda^k y^- + \mathcal{O}(\lambda^{2k})$. 
And last, we perform the \emph{scaling}
\[
x= -\frac{b}{d} \lambda^k X, \qquad y=-\frac{1}{d} \lambda^k Y, 
\]
under which the previous system becomes
\begin{equation}
\label{eq:gencase:final}
\begin{array}{rcl}
\bar{X} &=& Y + \mathcal{O}(\lambda^k) \\
\bar{Y} &=& M_1 + M_2 X - Y^2 \mathcal{O}(\lambda^k),
\end{array}
\end{equation}
with
$M_1=-d \lambda^{-2k} m_2 = -d \lambda^{-2k} \left( \mu + (cx^+ - y^-) \lambda^k + \mathcal{O}(k\lambda^k) \right)$ and $M_2=bc$, as it was claimed.

\qed

Lemma~\ref{lm:rescMar} shows that the limit form (that is, for large enough values of $k$ or, in other words, for close-enough orbits to $\Wsloc(O)$) for the first return map $T_{1k}=T_1 T_0^k$ (and similarly for $T_{2m}$) is the standard H\'enon map $\mathcal{H}$:
\[
\bar x = y, \qquad \bar y = M_1 + M_2 x - y^2,
\]
with Jacobian $J=-M_2 = -bc$. Recall that by~(\ref{detT12*}) and condition \textsf{[C]} we have $0<J<1$.
Bifurcations of fixed points of the standard H\'enon map are well known.
In the $(M_1, M_2)$-parameter plane, there are two bifurcation curves, namely
\begin{eqnarray*}
L^{+1}&:=\left\{ (M_1, M_2): 4 M_1 = - (1+ M_2)^2 \right\}, \\
L^{-1}&:=\left\{(M_1,M_2): 4 M_1 = 3 (1+M_2)^2 \right\},
\end{eqnarray*}
corresponding to the existence of a fixed point with a multiplier $+1$ (saddle-node fixed point) and a fixed point with a multiplier $-1$
(period doubling bifurcation), respectively. For $-1< M_2< 0$, the H\'enon map has no fixed points below the curve $L^{+1}$, has a stable (sink) fixed point in
the region between the bifurcation curves $L^{+1}$ and $L^{-1}$, while at $L^{-1}$ a period doubling bifurcation takes place and a stable 2-periodic orbit
appears above the curve $L^{-1}$.

Thus, using the relation~(\ref{eq:M1M2}) between the rescaled and the initial parameters we find that
\begin{eqnarray*}
\mu_k^+ &= \lambda^k (c \alpha  + \beta + \rho_k) + \frac{(1-bc)^2}{4d} \lambda^{2k}, \\[1.2ex]
\mu_k^- &= \lambda^k (c \alpha  + \beta + \rho_k)  - \frac{3(1-bc)^2}{4d} \lambda^{2k},
\end{eqnarray*}
where $\rho_k = \mathcal{O}(k\lambda^k)$ is small, $\alpha,\beta$ have been defined in~(\ref{xeqy}) and $b,c,d$ are Taylor coefficients of the map $T_1$ (see~(\ref{T12})).
This completes the proof of Theorem~\ref{th:th1}.

\qed

\begin{remark}
\label{rmk:non-intersec}
In general, the intervals $\delta_k$ do not intersect each other for different sufficiently large $k$.
However, when $c \alpha + \beta=0$, they can intersect and even appear nested. In the latter case, this implies that the diffeomorphism $f_0$ can possess simultaneously infinitely many periodic sinks and sources of all successive periods beginning from some (sufficiently) large number. This is a more delicate problem and it is out of the scope of this work. We recall that such phenomenon of ``global resonance'' with elliptic points was introduced in~\cite{GS01}
for area-preserving maps with homoclinic tangencies (see also \cite{GG09,
DGG15}).
\end{remark}

\section{Proof of Theorem~\ref{th:th2}}
\label{sec:proofTh2}

This proof will follow similar ideas and techniques as those employed in the proof of Theorem~\ref{th:th1}. We begin
by taking on $U_0$ the local $C^{r-1}$-coordinates $(x,y)$ provided by Lemma~\ref{LmSaddle}.
Recall that in these local coordinates the homoclinic points are denoted by $M_1^+=(x_1^+,0)$, $M_1^-=(0,y_1^-)$ in $\Gamma_1$ and $M_2^+=(x_2^+,0)$ and $M_2^-=(0,y_2^-)$ in $\Gamma_2$. They satisfy that $L(M_1^+) = M_2^-$ and
$L(M_1^-)=M_2^+$ (locally) since $R(M_1^+) = M_2^-$, $R(M_1^-)=M_2^+$, respectively. Now we consider the first return map $\Totkm= T_2 T_0^m T_1 T_0^k$ for single-round periodic 12-orbits. Thus, the following result holds:

\begin{lm}
%\label{lm:rescMar:th2}
Let us consider the family $\left\{f_\mu \right\}$ of Theorem~\ref{th:th2}, satisfying conditions~\textrm{\textsf{[A,B,C]}}.
Then, for large enough values of $k,m$, with $k \simeq m$, the %(figure-8)
first return map $\Totkm:
\sigma_k^0 \rightarrow \sigma_k^0$ can be brought, by a linear change of coordinates and
a suitable rescaling, to a reversible map asymptotically close as $k,m\to\infty$ to an area-preserving (symplectic) map of the form (see also \cite{DGGLS13}):
\begin{equation}
H\;:\left\{
\begin{array}{rl}
\bar x &= \Mt + \tilde c x -  y^2, \\
\tilde c \bar y &= -\Mt + y + \bar x^2,
\end{array}
\right.
\label{frmex}
\end{equation}
where
\begin{equation}
\tilde c = \frac{c}{b}\lambda^{k-m},\qquad
\Mt  = - \frac{d}{b^2}\lambda^{-2m}\left(\mu + c \lambda^k \beta + \lambda^m \alpha + O(k\lambda^k + m\lambda^m) \right).
\label{resc12:th2}
\end{equation}
The constants $\alpha,\beta$ are defined in~(\ref{xeqy}) and $b,c,d$ in expression~(\ref{T12}).
\end{lm}

From hypotheses~\textsf{[A]} and~\textsf{[C]} it follows that $\lambda >0$ and also
$\tilde c <0$ in the orientable case (if $T_1$ is orientable) and $\tilde c >0$ in
the non-orientable case (if $T_1$ is non-orientable).

\medskip

\begin{proof}
First, let us remind how coordinates are denoted on each domain around the homoclinic points $M_{1,2}^-$.
Thus, $(x,y)$-coordinates on $\Pi_i^{+}$ are denoted by $(x_{0i},y_{0i})$ and  by
$(x_{1i},y_{1i})$ on $\Pi_i^{-}$, for $i=1,2$. From Lemma~\ref{LmLocalMap}, the map $T_{0}^k: \Pi_2^+ \to \Pi_1^-$
will be defined on the strip $\sigma_k^{021}\subset\Pi_2^+$ and $T_0^k(\sigma_k^{021}) =
\sigma_k^{121}\subset\Pi_1^-$. Analogously, there exist strips  $\sigma_k^{011},
\sigma_k^{012}\subset\Pi_1^+$, and $\sigma_k^{022}\subset\Pi_2^+$ such that
$T_0^k(\sigma_k^{011}) = \sigma_k^{111}\subset\Pi_1^-$, $T_0^k(\sigma_k^{012}) =
\sigma_k^{112}\subset\Pi_2^-$ and $T_0^k(\sigma_k^{022}) = \sigma_k^{122}\subset\Pi_2^-$ (see
Fig.~\ref{fig:locmaps} for a comprehensive plot).
The first return map $\Totkm$ is given by the following chain of compositions:
\[
\sigma_k^{021}\stackrel{T_{0}^k}{\longmapsto}\sigma_k^{121}\stackrel{T_1}{\longmapsto}\sigma_m^{012}
\stackrel{T_{0}^m}{\longmapsto}\sigma_m^{112}\stackrel{T_2}{\longmapsto}\sigma_k^{021}
\]
(for a geometrical illustration see Fig.~\ref{fig:frtms(rev)2}).
These relations can be expressed in coordinates through the following set of equations ($T_0^k$, $T_1$, $T_0^m$, and $T_2$, respectively):
\begin{equation}
\begin{array}{rl}
x_{11} =& \lambda^k x_{02} (1 + k\lambda^k h_k(x_{02},y_{11}) ) \\
y_{02} =& \lambda^k y_{11} (1 + k\lambda_1^k h_k(y_{11},x_{02}) ), \\ \\
x_{01} - x_1^+ =& F_{1}(x_{11},y_{11}-y_1^-,\mu)\equiv \\
               &a x_{11} + b (y_{11}-y_1^-) +  \varphi_1(x_{11},y_{11},\mu),\\
y_{01} = & G_{1}(x_{11},y_{11}-y_1^-,\mu) \equiv\\
         &\mu + c x_{11} + d (y_{11}-y_1^-)^2  + \varphi_2(x_{11},y_{11},\mu), \\ \\
x_{12} =& \lambda^m x_{01} (1 + m\lambda^m h_m(x_{01},y_{12}) )\\
y_{01} =& \lambda^m y_{12} (1 + m\lambda^m h_m(y_{12},x_{01}) ), \\ \\
x_{12} =& G_{1}(\bar y_{02},\bar x_{02}-x_2^+,\mu) =  \\
        &\mu + c \bar y_{02} + d (\bar x_{02}-x_2^+)^2 +  \varphi_2(\bar y_{02},\bar x_{02},\mu),\\
y_{12} - y_2^- =& F_{1}(\bar y_{02},\bar x_{02}-x_2^+,\mu) =\\
                &a \bar y_{02} + b (\bar x_{02}-x_2^+) + \varphi_1(\bar y_{02},\bar x_{02},\mu).
\end{array}
\label{T12k}
\end{equation}
Observe that these formulas are presented in two different forms. Indeed, the local maps $T_0^{k,m}$ are given in cross-form while the global maps $T_{1,2}$ are written in explicit form. Thus, our first-return map $\Totkm$ can be defined, in cross-variables, as
$\Totkm: (x_{02},y_{11}) \mapsto (\bar{x}_{02},\bar{y}_{11})$, through the equation $\bar{y}_{02} = \lambda^k \bar y_{11} (1 +
k\lambda_1^k h_k(\bar y_{11},\bar x_{02}) )$ which plays an intermediate r\^ole.
As we did in Lemma~\ref{lm:rescMar}, we introduce new variables
\[
x_1 = x_{01} - x_1^+, \quad x_2 = x_{02} - x_2^+, \quad y_1 = y_{11}-y_1^-, \quad y_2
= y_{12}-y_2^-
\]
and rewrite system~(\ref{T12k}) as follows:
\begin{equation}
\begin{array}{l}
x_{1} = b y_{1} + \mathcal{O}(\lambda^k) + \mathcal{O}(y_1^2), \\ \\
\lambda^m y_{2}(1+m\lambda^{2m}\mathcal{O}(|x_1| + |y_2|)) = \\
\qquad (\mu + c \lambda^k x_{2}^+ - \lambda^m y_{2}^- + \mathcal{O}(k\lambda^{2k}+m\lambda^{2m})) + c \lambda^k x_2 +  d y_{1}^2  + \\
\qquad \mathcal{O}(\lambda^{2k}|x_2| + \lambda^k |x_2y_1| + |y_1|^3), \\ \\
\lambda^m x_{1}(1+m\lambda^{2m}\mathcal{O}(|x_1| + |y_2|)) =  \\
\qquad (\mu + c \lambda^k y_{1}^- - \lambda^m x_1^+ + \mathcal{O}(k\lambda^{2k}+m\lambda^{2m})) + c \lambda^k \bar y_1 + d \bar x_2^2 + \\
\qquad \mathcal{O}(\lambda^{2k}|\bar x_2| + \lambda^k |\bar x_2\bar y_1| + |\bar y_1|^3)\,\\ \\
y_{2} = b \bar x_{2} + \mathcal{O}(\lambda^k) + O(\bar x_2^2),
\end{array}
\label{T12k(2)}
\end{equation}
Take $x_1$ and $y_2$ from the first and fourth equations of (\ref{T12k(2)}) and substitute them in
the second and third ones. After this, we obtain the map $\Totkm: (x_2,y_1)\mapsto
(\bar x_{2},\bar y_{1})$ given in the following implicit form
\begin{align*}
%\label{T12k(3)}
%\lambda^m b \bar x_2 (1 + m\lambda^{m} \mathcal{O}(\bar x_2)) =& M  + d y_{1}^2  + c \lambda^k x_2 + \nonumber \\
%&\mathcal{O}(\lambda^{2k}|x_2| + \lambda^k |x_2y_1| + |y_1|^3), \\
&\lambda^m b \bar x_2 (1 + m\lambda^{m} \mathcal{O}(\bar x_2)) = \\
&\qquad M  + d y_{1}^2  + c \lambda^k x_2 + \mathcal{O}(\lambda^{2k}|x_2| + \lambda^k |x_2y_1| + |y_1|^3), \\
&\lambda^m b y_{1}(1 + m\lambda^{m} \mathcal{O}(y_1)) = \\
&\qquad M + c \lambda^k \bar y_1 + d \bar x_2^2 +  \mathcal{O}(\lambda^{2k}|\bar x_2| + \lambda^k |\bar x_2\bar y_1| + |\bar y_1|^3), \nonumber
\end{align*}
where
${\ds M =   \mu + c \lambda^k y_1^- - \lambda^m x_1^+ + O(k\lambda^{2k} + m\lambda^{2m})}$ or, equivalently,
\[
M =   \mu + c \lambda^k \beta + \lambda^m \alpha + O(k\lambda^{2k} + m\lambda^{2m}).
\]
Take into account that $x_1^+=-\alpha<0$ and $y_1^-=\beta$ (see formulas~(\ref{xeqy}) have been used.
Notice that up to this point, the procedure is symmetric. That is, we could have started our first-return map with $T_0^m$ instead of $T_0^k$ and the formulas would have been the same. This is
reflected in the fact that all the equations up to now, including the definition of the constant $M$, are invariant under $k \leftrightarrow m$.
Following the same procedure performed in the proof of Theorem~\ref{th:th1}, we rescale the coordinates. Indeed, consider
\[
x_2 = -\frac{b}{d}\lambda^m x, \qquad  y_1 = -\frac{b}{d}\lambda^m y,
\]
which bring the first return map $T_{12k}$ into the following rescaled form
\[
\begin{array}{l}
\bar x = \Mt + \tilde c x - y^2 + O(\lambda^k + \lambda^{2k-m}), \\
y = \Mt + \tilde c \bar y - \bar x^2 + O(\lambda^k + \lambda^{2k-m}),
\end{array}
\]
where $\tilde c$ and $\Mt$ satisfy (\ref{resc12:th2}).
%\begin{equation}
%\tilde c = \frac{c}{b}\lambda^{k-m}, \quad \Mt  = - \frac{d}{b^2}\lambda^{-2m}\left(\mu + c
%\lambda^k \beta + \lambda^m \alpha + O(k\lambda^{2k} + m\lambda^{2m}) \right),
%\label{resc12}
%\end{equation}
%as claimed.
This ends the proof of the lemma.
\end{proof}

\qed

\bigskip

\begin{figure}[bht]
\begin{center}
\includegraphics[width=14cm]
%,height=7cm]
{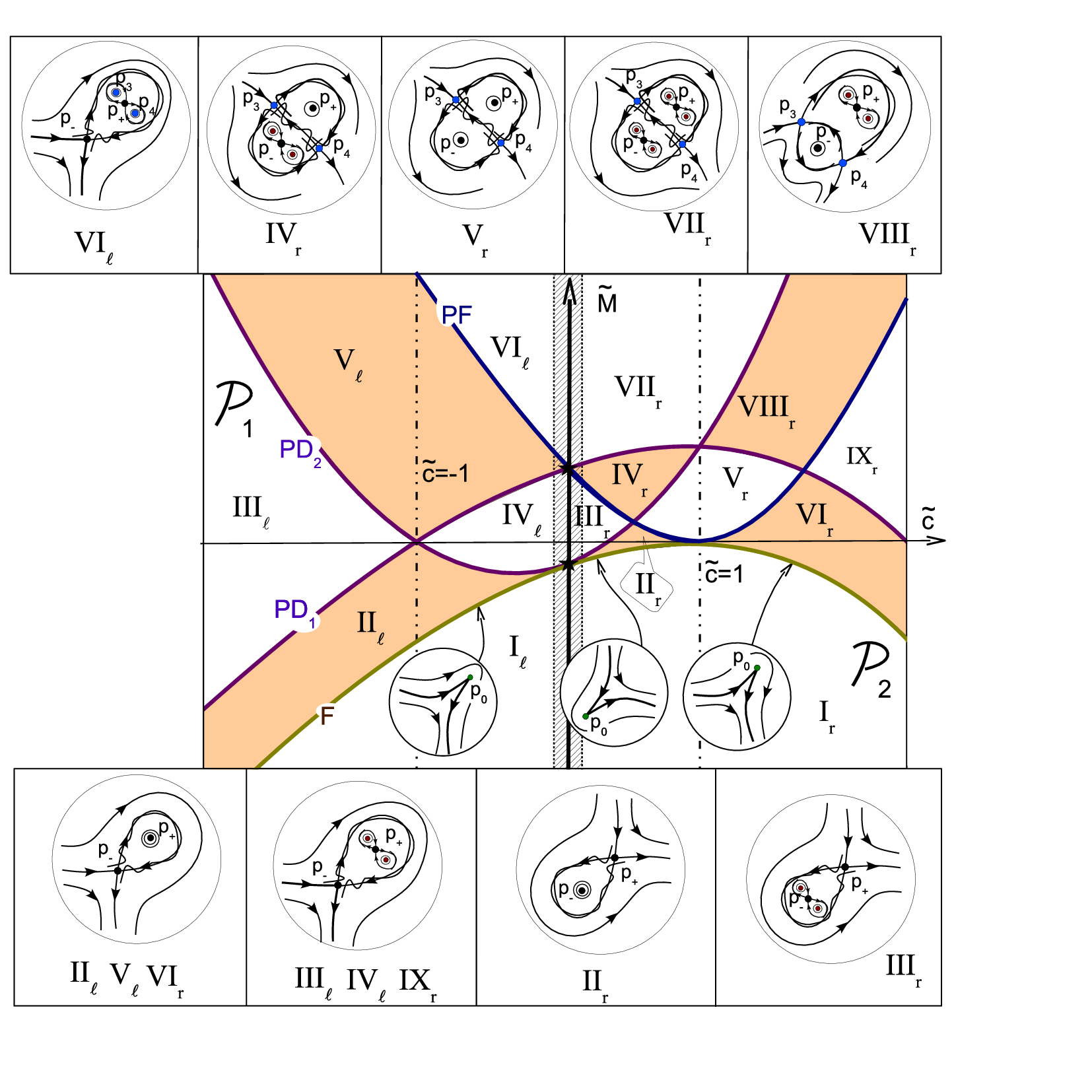} \caption{Elements of the bifurcation diagram for the map $H$: painted
regions correspond to the existence of symmetric elliptic and saddle fixed points of $H$. }
\label{Figdiag31}
\end{center}
\end{figure}

To complete the proof of Theorem 2
we need to detect the bifurcation boundaries of the intervals $\delta_{km}^c$. Since at $\mu\in\delta_{km}^c$ the first return map $\Totkm$ has two symmetric fixed points, one elliptic and another saddle, such boundaries can be found from the corresponding analysis of the map~(\ref{frmex}).
The bifurcation diagram for the symmetric fixed points of map~(\ref{frmex}) is shown in Fig.~\ref{Figdiag31}. We notice that it is essentially as the one in~\cite[page 16]{DGGLS13}. However, for the goals of \cite{DGGLS13}, 
searching only for symmetric fixed points was not sufficient, since the main problem there was to study symmetry breaking bifurcations (leading to the birth of a symmetric couple sink-source fixed points). 
This is not necessary here because the symmetric breaking bifurcations have been already determined in Theorem 1.

Like in~\cite[page 16]{DGGLS13}, the equations of the bifurcation curves $F$ (symmetric fold bifurcation), $P\!D_1$ and $P\!D_2$ (symmetric period doubling) and $P\!F$ (symmetry breaking pitch-fork) are
the following:
\begin{equation}
\begin{array}{rl}
F_0:  &   {\ds \Mt = -\frac{1}{4}\left(\tilde c -1 \right)^2},   \\\\
P\!D_1: &   {\ds \Mt = 1 -\frac{1}{4}\left(\tilde c -1 \right)^2},   \\\\
P\!D_2: &   {\ds \Mt = \frac{(\tilde c +1)(3\tilde c -1)}{4}},   \\\\
P\!F:   &   {\ds \Mt = \frac{3}{4}\left(\tilde c -1 \right)^2}.    \\\\
%PD(\textrm{asym}):   &  {\ds \Mt = \frac{(1-3\tilde c)(3 - \tilde c)}{4}}.
\end{array}
\label{bcurorient}
\end{equation}
These curves have the same equations for the orientable case, corresponding to the half-plane $\mathcal{P}_1 =\{\tilde c < -\varepsilon\}$ of the $(\tilde c,\Mt)$-parameter plane, and for the non-orientable case, corresponding to the half-plane $\mathcal{P}_2 =\{\tilde c > \varepsilon\}$, with an arbitrary small $\varepsilon>0$. Note that if $\tilde c =0$, then $c=0$ and therefore $T_1$ is not a diffeomorphism. So we exclude from the analysis a thin strip along the axes $\tilde c =0$ (the dashed strip in Fig.~\ref{Figdiag31}).

The curves (\ref{bcurorient}) divide the half-plane $\mathcal{P}_1$ in 6 domains $I_\ell,\ldots,V\!I_\ell$ and the half-plane $\mathcal{P}_2$ in 9 domains $I_r,\ldots,I\!X_r$. From these domains, we select two domains $I\!I_\ell$ and $V_\ell$ belonging to $\mathcal{P}_1$ and four domains $I\!I_r$, $I\!V_r$, $V\!I_r$ and $V\!I\!I\!I_r$ belonging to $\mathcal{P}_2$ which correspond to those values of the rescaled parameters $(\tilde c,\Mt)$ at which the map $H$ (and also the corresponding first return map $\Totkm$) has two symmetric fixed points: one saddle and another elliptic.
Note that for a given map $\Totkm$ the value of the parameter $\tilde c$ is uniquely determined.
%(see Remark~\ref{rem:invc}).

Then, the interval $\delta_{km}^c$ of values of the parameter $\mu$ corresponds to one of the intervals $\Delta_{\tilde c}| \tilde c = \mbox{const}$ of values of $\Mt$ that intersects some of the selected domains from its lower to its upper boundaries.

For instance, let us compute in the orientable case ($\tilde c <0$) the corresponding intervals $\delta_{km}^c$ of values of  $\Mt$ for the domain $I\!I_\ell$: 
\begin{equation}
\label{eq:IIelle}
\delta_{km}^c = \left(-\frac{1}{4}(\tilde c-1)^2, 1 -\frac{1}{4}(\tilde c-1)^2\right)\qquad\mbox{for}\;\;\tilde c \leq -1,
\end{equation}
and
\[
\delta_{km}^c = \left(-\frac{1}{4}(\tilde c-1)^2, \frac{1}{4}(\tilde c +1)(3\tilde c -1)\right)\qquad\mbox{for}\;\;-1<\tilde c <-\varepsilon.
\]
In both cases, the lower boundary corresponds to the symmetric fold bifurcation and the upper one to the symmetric period doubling.

Analogously, let us compute in the  non-orientable case ($\tilde c >\varepsilon$) the corresponding intervals $\delta_{km}^c$ for the domains $I\!I_r$ and $V\!I_r$: 
\[
\begin{array}{ll}
\delta_{km}^c = \left(-\frac{1}{4}(\tilde c-1)^2, \frac{1}{4}(\tilde c +1)(3\tilde c -1)\right) & \mbox{for}\;\;\tilde \varepsilon<c \leq 1/2; \\\\
\delta_{km}^c = \left(-\frac{1}{4}(\tilde c-1)^2, \frac{3}{4}(\tilde c -1)^2\right)  &\mbox{for}\;\; 1/2 < \tilde c  < 2 \;\;\mbox{and}\;\; \tilde c \neq 1; \\\\
\delta_{km}^c = \left(-\frac{1}{4}(\tilde c-1)^2, 1 - \frac{1}{4}(\tilde c -1)^2\right) &\mbox{for}\;\;  \tilde c  \geq  2 .
\end{array}
\]
In all three cases, the lower boundary corresponds to a symmetric fold bifurcation. However, the upper boundary corresponds to a symmetric period doubling for the first and the third case and to a symmetry breaking pitch-fork bifurcation for the intervals in the second case. 

We clearly will skip values of $k$ and $m$ such that $\tilde c =1$, that is, $\frac{c}{b}\lambda^{k-m}=1$. This is equivalent to say that $k - m = \frac{1}{\ln\lambda}\ln\frac{b}{c}$. 
Finally, we represent the intervals $\delta_{km}^c$ as intervals of values of $\mu$ using the relations (\ref{resc12:th2}). For example, for the intervals $\delta_{km}^c$ with $\tilde c  \leq -1$ (see~\eqref{eq:IIelle}), we obtain the following expressions for their bifurcation boundaries $\mu_{km}^{c+}\in F$ and $\mu_{km}^{c-}\in P\!D_1$:
\begin{eqnarray*}
\mu_{km}^{c+} =&& - c\lambda^k \beta - \lambda^m \alpha + \frac{b^2}{4d}(\tilde c -1)^2\lambda^{2m} \\
\mu_{km}^{c-} =&& - c\lambda^k \beta - \lambda^m \alpha + \frac{b^2}{d}\left(1 -\frac{1}{4}(\tilde c -1)^2\right)\lambda^{2m},
\end{eqnarray*}
and so on. Analogous explicit formulas can be obtained for the rest of the cases.

\qed

\section{Proof of Theorems 3 and 4
}
\label{sec:proofMainThm}

\subsection{Proof of Theorem 3}
\label{sec:sec6_1}

\begin{figure}[htb]
\begin{center}
\includegraphics[width=14cm]{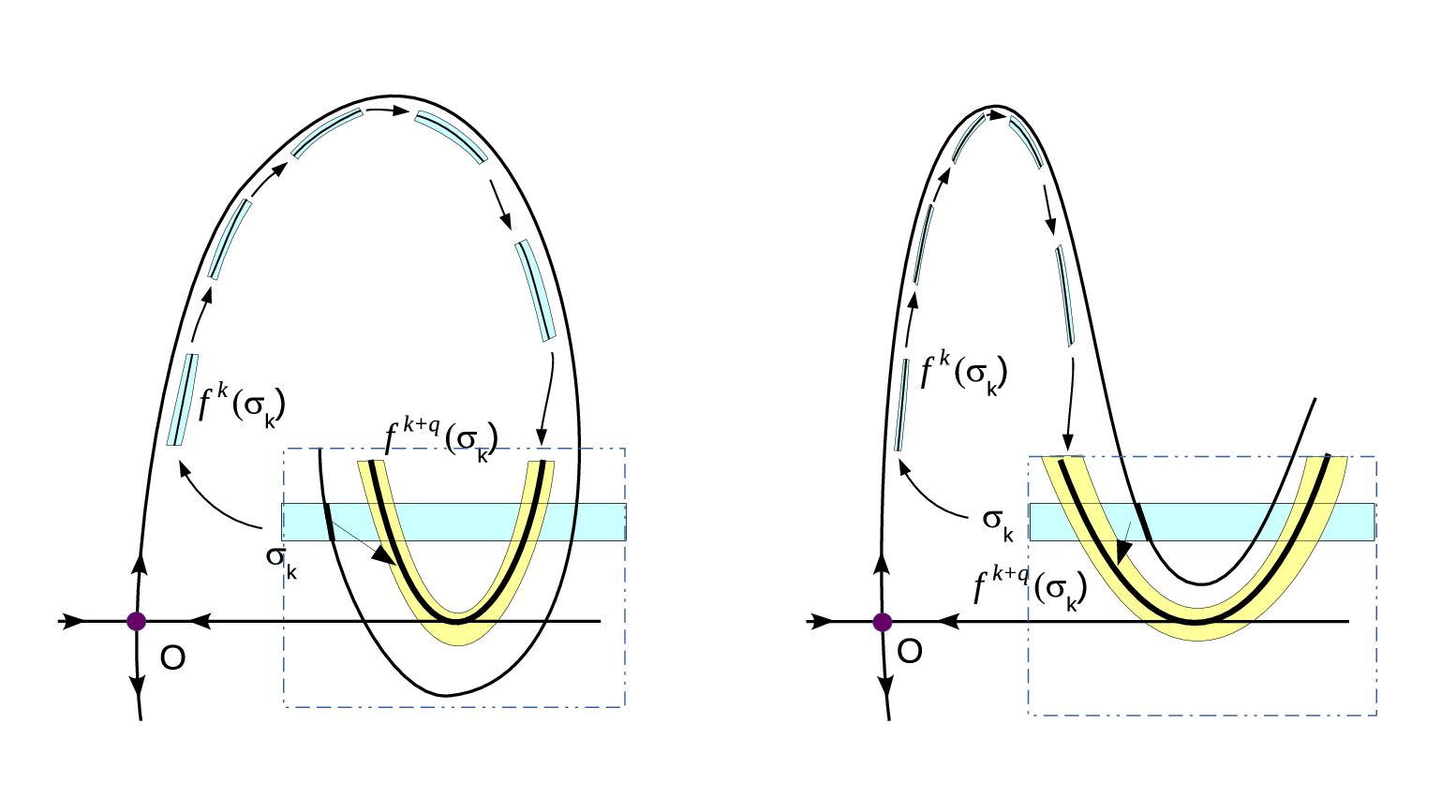}
\caption{Two examples of creation of secondary homoclinic tangencies to the point $O$ together with their Smale horseshoes}
\label{fig1newh}
\end{center}
\end{figure}

Its proof is quite standard (see, for instance,~\cite{N79,GST93b,GST99}). Namely, consider a single orbit $\Gamma_1$ and its neighbourhood $U_1$. From \cite{GST99} it is known that there exist $\left\{ \mu_k\right\}_k$, satisfying $\mu_k\to 0$ as $k\to\infty$, such that the map $f_{\mu_k}$ presents in $U_1$ a hyperbolic invariant set $\Lambda_k$ (a Smale horseshoe) such that (i) $W^u(\Lambda_k)$ is quadratically tangent to $W^s(O_\mu)$ and (ii) simultaneously, $W^u(O_{\mu_k})$ intersects transversally with $W^s(O_{\mu_k})$ (see Fig.~\ref{fig1newh}). Since all periodic points in $\Lambda_k$ have Jacobian less than 1 (by condition \textsf{[C]}) and, by the $\lambda$-Lemma, their stable and unstable manifolds accumulate (in a $C^r$-sense) to $W^s(O_\mu)$ and  $W^u(O_\mu)$, it follows that $\Lambda_k$ is a wild hyperbolic set (see~\cite{N79}). The latter assertion implies that, arbitrary close to $\mu=0$, there exist intervals of values of $\mu$ for which $W^u(\Lambda_k)$ and $W^s(\Lambda_k)$ have points with quadratic tangency.
Thus, one obtains that the values of $\mu$ for which the map $f_\mu$ has a nontransversal homoclinic orbit $\Gamma_{1\mu}\subset U_1$ are dense in these intervals.

\qed

\subsection{Proof of Theorem 4}
\label{sec:sec6_2}

The proof of this theorem follows from Theorems~\ref{th:th1} and~\ref{th:th2} and a standard procedure of embedding intervals applied to any arbitrary point belonging to any interval $n_i$ from Theorem 3. Indeed, take any $\bar\mu\in n_i$. Arbitrary close to $\bar\mu$ there is $\bar\mu_1\in n_i$ such that $f_{\bar\mu_1}$ has a couple of homoclinic tangencies of the initial type. Hence, by Theorem~\ref{th:th1}, near $\bar\mu_1$ there exists an interval $I_1 \subset n_i$ such that at $\mu \in I_1$ the diffeomorphism $f_\mu$ has a periodic couple ``sink-source''. In turn, since $n_i$ is the Newhouse interval, in $I_1$ we find an interval $I_2$ such that the diffeomorphism $f_\mu$ at $\mu\in I_2$ has simultaneously, a periodic couple ``sink-source'' (as $\mu\in I_1$) and a symmetric elliptic periodic orbit. Repeating this procedure beginning from the interval $I_2$ we obtain a sequence $I_2,I_4,...$ of embedding intervals such that at $\mu\in I_{2j}$ the diffeomorphism $f_\mu$ has  $j$ periodic couples ``sink-source'' and $j$ symmetric elliptic periodic orbits, etc.

\qed

\section{Some examples} \label{sec:ex}

In this section we provide some simple examples of planar reversible maps undergoing a ``fish'' or figure-$8$ quadratic homoclinic tangency. They are Poincar\'e maps of periodically perturbed planar reversible differential systems. By construction hypotheses \textsf{[A,B]} will be straightforwardly satisfied. The fulfilment of condition \textsf{[C]} is expected by numerical checking because of the large freedom one has to produce many close variants of the periodic perturbations.
The \emph{basic} systems will be the well-known \emph{Duffing equation} and the \emph{Cubic potential} (the ``fish''), both Hamiltonian and reversible. A similar approach was performed by Duarte in~\cite{Duarte00}.

\subsection{Perturbed Duffing equation}

Let us consider the vertical Duffing equation
\begin{equation}
\label{eqn:vert_duffing}
\left\{
\begin{array}{rcll}
\dot{x} &=& y - y^3 & + \varepsilon f(x,y,t) \\
\dot{y} &=& x       & + \varepsilon g(x,y,t).
\end{array}
\right.
\end{equation}
For $\varepsilon \ne 0$, system~(\ref{eqn:vert_duffing}) is Hamiltonian, reversible (with respect to linear involutions, $R(x,y)=(x,-y)$ and $S(x,y)=(-x,y)$) and presents a couple of ($R$-)symmetric homoclinic solutions to the origin. These figure-$8$ homoclinic curves (single-round $12$-orbits)
can be parameterized by $\Gamma_h^-(t)=(x_h(t),\pm y_h(t))$, where
\begin{equation*}
%\label{def:duffing:homoclinic}
x_h(t)=-\sqrt{2} \sech(t) \tanh(t), \qquad
y_h(t)=\sqrt{2} \sech(t)
\end{equation*}
for $t \in (-\infty,+\infty)$.
Moreover, the following properties hold:
(i) $x_h(t)=\dot{y}_h(t)$; (ii) $(x_h(0),y_h(0))=(0,\sqrt{2})$; (iii) $y_h(t)$ has a pole of order $1$ at the points $\pm \pi \rmi /2$ (and, therefore, $x_h(t)$ has poles of order $2$ at the same points).

Our aim is to provide some examples of periodic perturbation of~(\ref{eqn:vert_duffing}), preserving $R$-reversibility and not in general the Hamiltonian character, such that the homoclinic invariant curves of the origin undergo a quadratic tangency (and, therefore, infinitely many of them).
It is straightforward to check that, for $\varepsilon \ne 0$, system~(\ref{eqn:vert_duffing}) is $R$-(time) reversible if and only if
$f(x,-y,-t)= -f(x,y,t)$ and $g(x,-y,-t)= g(x,y,t)$.
The existence of (tangent) quadratic homoclinic points will be carried out by selecting a
simple suitable perturbation and parameters $\omega_j, t_0^*$ such that the corresponding \emph{Melnikov function} $M(t_0)$ has a double-zero at $t_0=t_0^*$. Melnikov function is given by
\begin{equation*}
%\label{eqn:Melnikov:duffing}
M(t_0)= \int_{-\infty}^{+\infty} \left(F \wedge G \right) \left( x_h(t),y_h(t),t+t_0 \right) \, dt,
\end{equation*}
where
\[
F(x,y)= \left(
\begin{array}{c}
y-y^3 \\
x
\end{array}
\right), \qquad
G(x,y,t)= \left(
\begin{array}{c}
f(x,y,t) \\
g(x,y,t)
\end{array}
\right)
\]
and $F\wedge G = (y-y^3) g(x,y,t) - x f(x,y,t)$. To produce such example, we restrict ourselves to the case where $g\equiv 0$ and $f(x,y,t)$ a (periodic) linear combinations
of \emph{odd} functions of the form $x \sin \omega t$, that is,
\begin{equation*}
\left\{
\begin{array}{rcll}
\dot{x} &=& y - y^3 & + \varepsilon x \sum_{j=0}^N a_j \sin \omega_j t\\
\dot{y} &=& x       &
\end{array}
\right.
\end{equation*}
with commensurable $\omega_0, \omega_1, \ldots, \omega_N$. Having in mind that $x_h^2(t) \sin \omega_j t$ is an odd function in $t$ (and, therefore, its integral
over $(-\infty,+\infty)$ is null) it follows that
\begin{eqnarray*}
\lefteqn{ M(t_0) = - \sum_{j=0}^N a_j \int_{-\infty}^{+\infty}  x_h^2(t) \sin \omega_j(t+t_0) \, dt = } \\
&{\ds  - \sum_{j=0}^N a_j \left( \int_{-\infty}^{+\infty}  x_h^2(t) \cos \omega_j t \, dt \right) \sin \omega_j t_0 = }\\
&{\ds  - \frac{\rme^{\pi/2}}{3 \sinh(\pi/2)} \sum_{j=0}^N \left( a_j \sinh\left( \frac{\pi \omega_j}{2} \right) \, (\omega_j^2 - 2) \omega_j \sin \omega_j t_0 \right),}
\end{eqnarray*}
provided by the \emph{residues} integration
\[
\mathrm{\,Res}\left( x_h^2(t) \cos \omega_j t,t= \frac{\pi \rmi}{2} \right) =
\frac{\rme^{\pi/2}}{3 \sinh(\pi/2)} \sinh\left( \frac{\pi \omega_j}{2} \right) \, (\omega_j^2 - 2) \omega_j.
\]
Let us consider as a particular example, the case $\omega_0=1, \omega_1 = \omega \in \Z \setminus \left\{ 1 \right\}$, $a_0=\alpha$ and $a_1=\beta$ with $\alpha\beta \ne 0$.
Indeed,
\begin{equation*}
\left\{
\begin{array}{rcll}
\dot{x} &=& y - y^3 & + \varepsilon x \left( \alpha \sin t + \beta \sin \omega t \right)\\
\dot{y} &=& x.       &
\end{array}
\right.
\end{equation*}
Now, the Melnikov function reads
\[
M(t_0)=  - \frac{\rme^{\pi/2}}{3\sinh(\pi/2)} \left( - \alpha \sinh\left( \frac{\pi}{2} \right) \, \sin t_0  \right. \\[1.2ex]
 \left. + \beta \sinh \left( \frac{\pi\omega}{2} \right) \, (\omega^2-2) \omega \, \sin \omega t_0 \right).
\]
We seek for values of $\omega$ and $t_0$ satisfying that $M(t_0)=M'(t_0)=0$ and $M''(t_0)\ne 0$, i.e., giving rise to a quadratic homoclinic tangency.
Denoting $A=\alpha \sinh(\pi/2)$ and $B_{\omega}= \beta \sinh(\pi\omega/2) \, (\omega^2 -2) \omega $, this is equivalent to look for double zeroes of $\varphi_{\omega}(t_0)= -A \sin t_0 + B_{\omega} \sin \omega t_0$. Since $\beta\ne 0$ and $\omega \ne 0$ it turns out that $B_{\omega}$ does not vanish as well. It is straightforward to check that $\varphi_{\omega}(t_0)=\varphi_{\omega}'(t_0)=0$, $\varphi_{\omega}''(t_0)\ne 0$ reduces to find $\omega$ and $t_0$ with $\omega t_0 \ne k \pi$, for $k\in \Z$, satisfying $A \sin t_0 = B_{\omega} \sin \omega t_0$ and $A \cos t_0 = \omega B_{\omega} \cos \omega t_0$.

It is simple to prove that there is no solution $t_0$ for $\omega=2$. Indeed, $\omega t_0 \notin \pi \Z$ implies that $t_0 \ne k\pi/2$ for $k\in \Z$. Imposing the two other conditions leads us, first, to $A=2B_2 \cos t_0$ and, second, to $\sin t_0=0$, a contradiction with the fact that $t_0 \ne k\pi/2$.
If we choose $\omega=3$ and (for instance) $t_0=\pi/2$, that is
\begin{equation*}
\left\{
\begin{array}{rcll}
\dot{x} &=& y - y^3 & + \varepsilon x \left( \alpha \sin t + \beta \sin 3 t \right)\\
\dot{y} &=& x,       &
\end{array}
\right.
\end{equation*}
the latter conditions reduce to $B_3=-A$ and having in mind that
$A=\alpha \sinh(\pi/2)$ and $B_{3}= 21 \beta \sinh(3\pi/2)$ it follows that we have a quadratic homoclinic point at $t_0=\pi/2$ for $\omega=3$ provided
\[
\beta = - \frac{\sinh(\pi/2)}{21 \sinh(3\pi/2)} \, \alpha.
\]

\subsection{Perturbed ``fish'' equation}

This example of single-round 1- and 2-orbits, based on the \emph{fish} equation, is given by
\begin{equation*}
%\label{system:fish}
\left\{
\begin{array}{rcll}
\dot{x} &=& y & + \varepsilon f(x,y,t)\\
\dot{y} &=& x - x^2 & + \varepsilon g(x,y,t).
\end{array}
\right.
\end{equation*}
For $\varepsilon=0$ this fish equation is (time) $R$-reversible, with $R$ the involution $(x,y)\mapsto (x,-y)$, and presents a ($R$)-symmetric homoclinic solution to the origin, namely, $\Gamma_h(t)=(x_h(t),y_h(t))$, where
\[
x_h(t)=\frac{\sqrt{3}}{2} \sech^2 \left( \frac{t}{2} \right), \quad
y_h(t)=\dot{x}_h(t)=-\frac{\sqrt{3}}{2} \sech^2 \left( \frac{t}{2} \right) \tanh \left( \frac{t}{2} \right)
\]
Function $x_h(t)$ has a pole of order $2$ at $\pm \pi \rmi$ and, therefore, $y_h(t)$ has them of order $3$.
If we ask the perturbation $(f,g)$ to preserve the $R$-reversibility, it must satisfy that
$f(x,-y,-t)=-f(x,y,t)$ and $g(x,-y,-t)=g(x,y,t)$.
Proceeding like in the previous example, the Melnikov function for a general reversible perturbation $(f,g)$ reads as follows
\begin{eqnarray*}
\lefteqn{M(t_0) = \int_{-\infty}^{+\infty} \left(F \wedge G \right) \left( x_h(t),y_h(t),t+t_0 \right) \, dt = } \\
&{\ds \int_{-\infty}^{+\infty}  y_h(t) g(x_h(t),y_h(t),t+t_0) \, dt  - \int_{-\infty}^{+\infty} (x_h(t)- x_h^2(t)) f(x_h(t),y_h(t),t+t_0) \, dt.}
\end{eqnarray*}
As before, we restrict ourselves to a simpler case, namely,
\[
f\equiv 0, \qquad \qquad  g(x,y,t)=g(y,t) = y \sum_{j=0}^N b_j \sin \omega_j t,
\]
again with $\omega_0, \omega_1, \ldots, \omega_N$ commensurables.
As we did for the Duffing equation, we select a simple example giving rise to a homoclinic quadratic point. Indeed, we choose $\omega_0=2$, $\omega_1=6$ (they are the smallest satisfying it), $t_0=\pi/4$ and denote $b_0=\alpha$, $b_1=\beta$ (with $\alpha \beta \ne 0$). Indeed,
\begin{equation*}
\left\{
\begin{array}{rcll}
\dot{x} &=& y & \\
\dot{y} &=& x - x^2 & + \varepsilon y \left( \alpha \sin 2t + \beta \sin 6t \right).
\end{array}
\right.
\end{equation*}
Thus, our Melnikov function reads
\begin{eqnarray*}
M(t_0)=& \frac{4}{5}\pi \left( \alpha \sinh(2\pi) \cdot (2^4-1)\cdot 2 \cdot \sin(2 t_0) \right.
\\ & \left.+ \beta \sinh(6\pi) \cdot (6^4-1) \cdot 6 \cdot \sin(6 t_0) \right),
\end{eqnarray*}
which can be written as $A \sin 2t_0 + B \sin 6t_0$ with
\[
A=\frac{4}{5}\pi \alpha \sinh(2\pi) \cdot (2^4-1) \cdot 2, \quad
B=\frac{4}{5}\pi \beta \sinh(6\pi) \cdot (6^4-1) \cdot 6.
\]
Taking $A=B$ it follows that $M(\pi/4)=M'(\pi/4)=0$ and $M''(\pi/4)=32B \ne 0$, which provides the condition
\[
\beta = \frac{(2^4-1) \sinh(2\pi)}{3 (6^4-1)\sinh(6\pi)}\,  \alpha.
\]

\section*{Acknowledgements}

The authors are grateful to D. Turaev and L. Lerman for fruitful discussions and useful comments. MG warmly thanks the Department of Mathematics 
of Uppsala University for their hospitality and support during her stay at Uppsala University. JTL thanks the Centre de Recerca Matem\`atica (CRM) for its hospitality.

\end{document}